\documentclass{amsart}
\newif\ifcolor
\colortrue

\usepackage{amssymb,amsmath,amsthm,epsfig,psfrag}
\ifcolor
  \usepackage[dvips]{color}
  \definecolor{dkgreen}{rgb}{0,0.45,0.08}
  \newcommand{\Red  }[1]{{\color{red}{#1}}}
  \newcommand{\Blue }[1]{{\color{blue}{#1}}}
  \newcommand{\DkGrn}[1]{{\color{dkgreen}{#1}}}

\else
  \newcommand{\Red  }[1]{#1}
  \newcommand{\Blue }[1]{#1}
  \newcommand{\DkGrn}[1]{#1}

\fi

\theoremstyle{plain}
\newtheorem{thm}[equation]{Theorem}
\newtheorem{prop}[equation]{Proposition}
\newtheorem{lemma}[equation]{Lemma}
\newtheorem{cor}[equation]{Corollary}

\theoremstyle{definition}
\newtheorem{remark}[equation]{Remark}
\newtheorem{example}[equation]{Example}
\newtheorem{defn}[equation]{Definition}

\newcommand{\claim}{\smallskip\noindent{\bf Claim:\ }}

\newcommand{\defterm}[1]{\emph{#1}}

\numberwithin{equation}{section}

\newcommand{\CC}{{\mathcal C}}
\newcommand{\MM}{{\mathcal M}}
\newcommand{\Ct}{C_{\theta}}

\renewcommand{\Re}{\mathop{\rm Re}\nolimits}
\renewcommand{\Im}{\mathop{\rm Im}\nolimits}

\newcommand{\R}{{\mathbb R}}
\newcommand{\Z}{{\mathbb Z}}

\newcommand{\En}{{\bf E}_n}
\newcommand{\On}{{\bf O}_n}
\newcommand{\Bodd}{{B^{\bf o}}}
\newcommand{\Beven}{{B^{\bf e}}}
\newcommand{\Modd}{{M^{\bf o}}}
\newcommand{\Meven}{{M^{\bf e}}}

\newcommand{\Match}{M}
\newcommand{\Basket}{B}
\newcommand{\st}{\;:~}

\newcommand{\BB}{{\mathcal{B}}}
\newcommand{\QQ}{{\mathcal{Q}}}

\begin{document}
\title{Harmonic algebraic curves and noncrossing partitions}
\author[Martin, Savitt and Singer]{Jeremy Martin, David Savitt, and Ted Singer}

\address{Jeremy L. Martin, Department of Mathematics, University of Kansas, 405 Snow Hall,
1460 Jayhawk Blvd., Lawrence, KS 66045}
\email{jmartin@math.ku.edu}
\address{David Savitt, Department of Mathematics, University of Arizona, 617 N.\
Santa Rita Ave., Tucson, AZ 85721}
\email{savitt@math.arizona.edu}
\address{Ted Singer, 38 Greenacres, Scarsdale, NY 10583}
\email{ted.singer@gmail.com}

\thanks{First author supported in part by an NSF Postdoctoral Fellowship.}

\begin{abstract}
Motivated by Gauss's first proof of the Fundamental Theorem
of Algebra, we study the topology of harmonic algebraic
curves.  By the maximum principle, a harmonic curve has no bounded components;
its topology is determined by the combinatorial data of a noncrossing matching.
Similarly, every complex polynomial gives rise to
a related combinatorial object that we call a \textit{basketball},
consisting of a pair of noncrossing matchings satisfying one
additional constraint.  
We prove that every noncrossing matching arises from some
harmonic curve, and deduce from this that every basketball
arises from some polynomial.
\end{abstract}

\maketitle

%===========================================================================

\section*{Introduction}

The first proof of the Fundamental Theorem of Algebra was
published in 1746 by d'Alembert \cite{dAlembert}.  In an attempt
to correct the lack of rigor in d'Alembert's approach, as well as
in subsequent attempts, Gauss offered a new proof in 1799 in his
doctoral thesis \cite{Gauss}. Gauss's argument, while
characteristically elegant, was itself not entirely satisfactory.
Both d'Alembert's and Gauss's proofs were subsequently made
completely rigorous (for the whole story, see, e.g.,
\cite[pp.~195--200]{Stillwell}).

\begin{figure}[h]
\begin{center}
\hfill\phantom{.}
\resizebox{1.5in}{1.5in}{\includegraphics{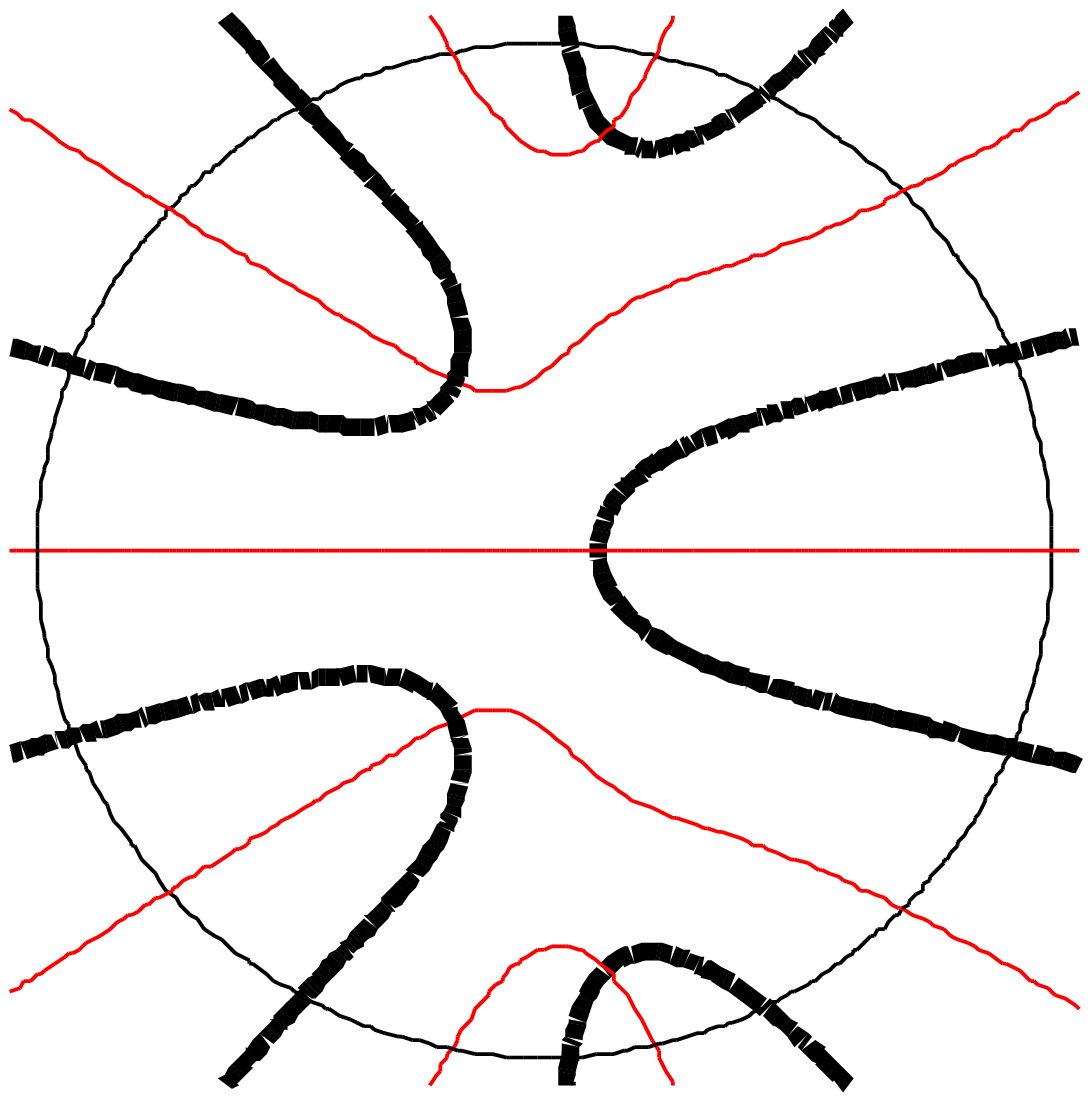}}
\hfill
\resizebox{1.5in}{1.5in}{\includegraphics{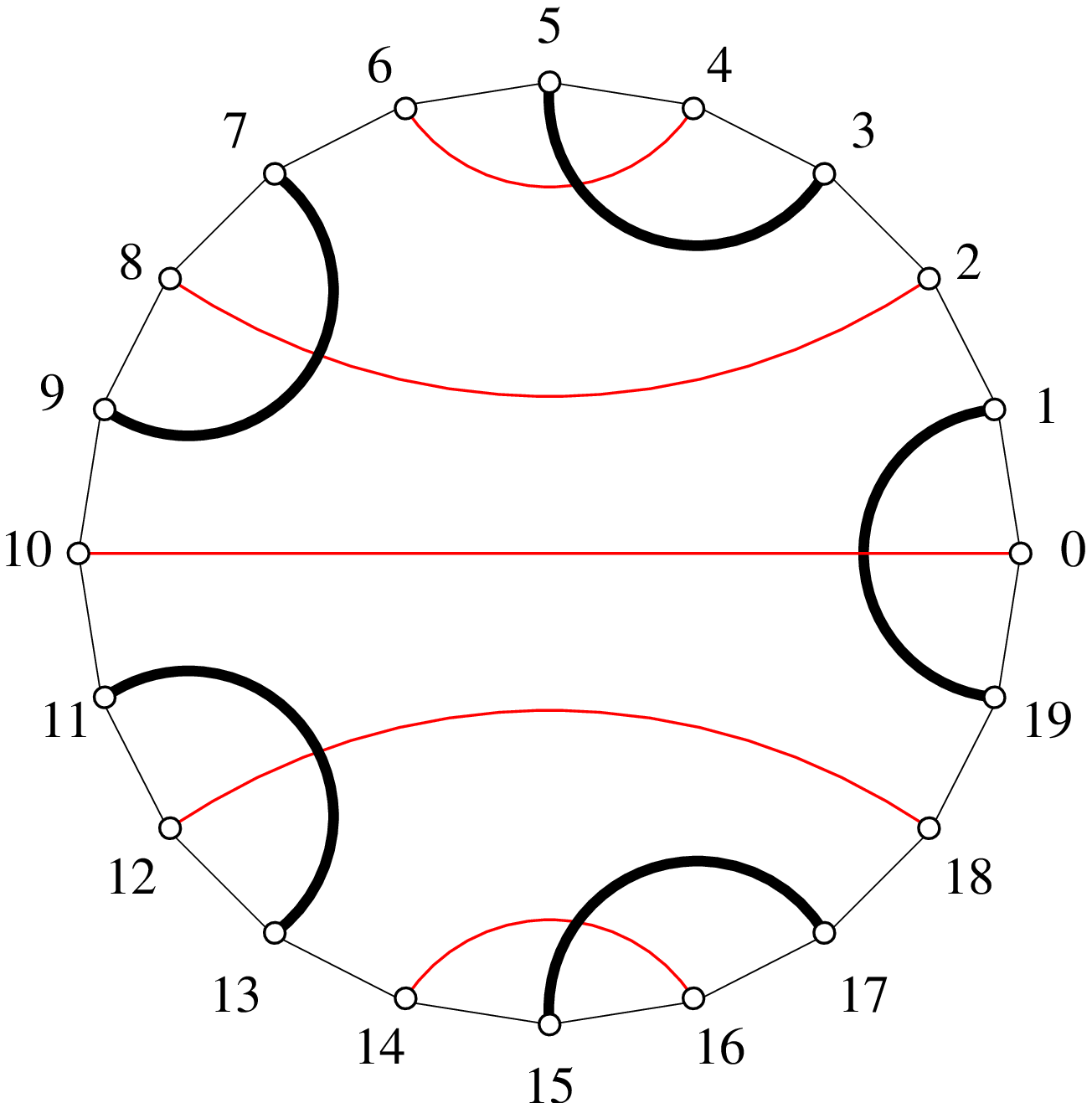}}
\hfill\phantom{.} %% this is a kludge to get LaTeX to space things out nicely
\end{center}
\caption{
  \emph{Left}: The curves $R$ ({\bf thick}) and $I$ (\Red{thin})
   arising from $f(z) = z^5+6z^3+3z^2+5z-2$, with a circle $S$ superimposed.
  \emph{Right}: The associated combinatorial basketball.
\label{basketball-example-figure}}
\end{figure}

Gauss approached the problem by examining the real algebraic curves
  \begin{align*}
    R &= \{z \st \Re f(z)=0\},\\
    I &= \{z \st \Im f(z)=0\},
  \end{align*}
in $\R^2$, where~$f(z)$ is a complex polynomial of degree~$n$.  (See Figure~\ref{basketball-example-figure} for an example.)  Note that
the set of roots of~$f$ is precisely $R\cap I$.  For a suitably
large disk $D$, each of the curves $R$ and $I$ meets the boundary circle $S=\partial
D$ in exactly $2n$ points, and it can easily be shown that the points of
$R\cap S$ alternate with those of $I\cap S$ around $S$.  Gauss argued
that $R\cap D$ (resp.\ $I\cap D$) must consist of $n$ components
and that each component of $R$ must
cross some component of $I$ inside $D$; therefore, $\#(R \cap I)\geq n$.

We are interested in the topology of the curves $R$ and $I$,
both singly and together.  Much of the research on
the topology of real plane algebraic curves (the topic of
Hilbert's sixteenth problem) has focused on classifying the possible
configurations of ovals (bounded connected components); see, e.g.,
\cite{Viro}.  However, the curves $R$ and $I$ are \emph{harmonic},
and so the maximum principle \cite[Theorem~21, p.~166]{Ahlfors} implies
that they have no ovals at all.  The topological information that can be extracted
from the pair $R,I$ is of a different sort entirely.

Suppose that $R$ and $I$ are both nonsingular.  (In particular, it is
necessary that $f(z)$ have no repeated roots; as we will see, this condition
is in a certain sense sufficient.)  Labeling the $2n$ points
of $R\cap S$ cyclically and pairing off those points that lie in
the same component of~$R$ yields a \emph{noncrossing matching} of
order $n$. (Noncrossing partitions are an important subject in modern combinatorics;
for an overview, see the excellent survey by Simion \cite{Simion}.)
Repeating the construction for $I$ instead of $R$ yields a second
noncrossing matching that is interlaced with the first.  The
matching obtained from $R$ must cross the matching obtained from
$I$ exactly $n$ times, corresponding to the $n$ roots of the
polynomial $f$.  Two arbitrary interlaced noncrossing matchings of
order $n$ must have at least $n$ crossings between them (Lemma
\ref{lemma:odd-cross}), and may have more; if the number of such
crossings is exactly $n$, we call this pair of noncrossing
matchings a \emph{basketball of order~$n$}.  For example, the right-hand side
of Figure~\ref{basketball-example-figure} is the basketball
  \begin{align*}
  \big( & \{\{0,10\},\{2,8\},\{4,6\},\{12,18\},\{14,16\}\},\\
  & \{\{1,19\},\{3,5\},\{7,9\},\{11,13\},\{15,17\}\}\big).
  \end{align*}

It is natural to ask whether every basketball arises from
a polynomial in this way.  Our main result, the Inverse
Basketball Theorem (Theorem~\ref{IBT}), answers this question
in the affirmative.  The proof is
constructive, and draws on elementary tools from combinatorics,
topology and complex analysis.

There is a natural bijection between basketballs of order $n$
and noncrossing partitions with $n$ blocks of size 4.  It
follows from a result of Edelman
\cite{Edelman} that the basketballs are enumerated by the
quasi-Catalan numbers $\frac{1}{3n+1}\binom{4n}{n}$.  
These numbers also count plane quaternary trees and dissections of a
polygon into pentagons; there are natural bijections between
basketballs and each of these combinatorial sets.
The quasi-Catalan numbers also enumerate certain plane trees (see
\cite{LW}), although there is no obvious bijection involved.

A more subtle combinatorial-topological invariant of the
polynomial $f$ arises from the curves
  $$\Ct(f) = \{ z \st \Im(e^{-i\theta} f(z)) = 0\},$$
regarded as a family parameterized by $\theta\in\R/\pi\Z$.
It turns out that $C_\theta(f)$ is singular for
only finitely many values of $\theta$.  Therefore, we can study
the family of noncrossing matchings, equipped with a cyclic order,
obtained by letting $\theta$ vary; we call this family the
\emph{necklace of matchings} associated to $f$.  Notice that the
necklace depends only on the polynomial $f$ itself, and not on
a choice of angle.

The paper is structured as follows.  In Section~\ref{section:FTA}, we give
a modern exposition of Gauss's proof of the Fundamental Theorem of Algebra,
and explain how Gauss's ideas may be used to associate basketballs
with complex polynomials.  Section~\ref{section:basketball}
is devoted to the combinatorics of basketballs.  Section~\ref{section:IBT}
contains the proof of the Inverse Basketball Theorem.
We conclude in Section 4 with some brief remarks on necklaces of
matchings; there appears to be much more to say here from both the combinatorial
and geometric points of view.

%===========================================================================

\section{Gauss's proof of the Fundamental Theorem of Algebra} \label{section:FTA}

We begin by describing Gauss's approach to the Fundamental Theorem
of Algebra \cite{Gauss}; for the technical details, see
Gersten and Stallings~\cite{Gersten}. Let $f(z)$ be a
monic polynomial of degree $n$, and consider the two curves
  \begin{align*}
    R &= \{z \st \Re f(z)=0\},\\
    I &= \{z \st \Im f(z)=0\} \,.
  \end{align*}
In polar coordinates $z=re^{i\theta}$, the curves $R$ and $I$ are
given by equations
  \begin{align*}
  \Re\,f(z) &= r^n \cos n\theta + \text{(lower terms)} = 0,\\
  \Im\,f(z) &= r^n \sin n\theta + \text{(lower terms)} = 0.
  \end{align*}
Let $S_r$ denote the circle $\{z\st|z|=r\}$, and $D_r$ the disk
$\{z\st|z|\le r\}$. By taking $r$ sufficiently large, we can
ensure that
  \begin{itemize}
  \item $\Re\,f(z)$ has $2n$ zeros on $S_r$, arbitrarily close (in angle) to
    the zeros of $\cos n\theta$; and
  \item $\Im\,f(z)$ has $2n$ zeros on $S_r$, arbitrarily close to
    the zeros of $\sin n\theta$.
  \end{itemize}

That is, $R \cap S_r$ and $I \cap S_r$ each contain exactly $2n$
points, and these $4n$ points alternate around the circle.
(See Figure~\ref{basketball-example-figure} for an example.)
Suppose that $R$ and $I$ are both smooth curves.  Then the
intersection $R\cap D_r$ must contain $n$ disjoint arcs
$R_1,\ldots,R_n$, each of which joins two distinct points of $R
\cap S_r$. Since these arcs do not cross, each arc $R_i$ must have
an even number of points of $R\cap S_r$ on either side of it,
by the Jordan Curve Theorem.

Likewise, $I\cap D_r$ contains $n$ disjoint arcs $I_1,\ldots,I_n$,
each of which joins two distinct points of $I \cap S_r$. Since
there is exactly one point of $I \cap S_r$ lying between each
point of $R \cap S_r$, each arc $R_i$ must have an \emph{odd}
number of points of $I \cap S_r$ on either side of it.  It follows
that some arc $I_j$ joins a point on one side of $R_i$ with a
point on the other, and therefore must intersect $R_i$.  These $n$
intersection points, one for each $R_i$, are the roots of our
polynomial $f$. (This step requires only the
Intermediate Value Theorem, rather than the Jordan Curve Theorem,
because $\Re f(z)$ changes sign along the arc $I_j$.)
This establishes:

\begin{prop} \label{prop:nonsing} Let $f(z)$ be a monic polynomial such that the curves
$R$ and $I$ are nonsingular.  Then $f(z)$ has $n$ roots.
\end{prop}

Gersten and Stallings \cite{Gersten} complete the proof of the
Fundamental Theorem of Algebra by proving that if $f(z)$ has no
roots, then $f(z)$ can be perturbed to obtain a new polynomial that still
has no roots but for which $R$ and $I$ are nonsingular.  We give a
different approach.

Observe that the curve $I$ has $2n$ asymptotes, with slopes at
angles $k\pi/n$ for $k=0,\ldots,2n-1$.  Similarly, $R$ has
asymptotes with slopes at angles $(k + \frac{1}{2})\pi/n$. More
generally, consider the family of real plane curves
  \begin{equation} \label{define-Ctheta}
  \Ct(f) = \{ z \st \Im(e^{-i\theta} f(z)) = 0\}\,,
  \end{equation}
parameterized by the circle $\R/\pi\Z$.  We will often abuse notation
by identifying $\theta\in\R/\pi\Z$ with its coset representative in the interval
$[0,\pi)$.  Note that $\Ct(f)$ has $2n$ asymptotes, with slopes at angles
$\frac{k\pi + \theta}{n}$ for $k=0,\ldots,2n-1$.  We regard
$\arg(z)$ as an element in $\R/2\pi\Z$.

\begin{figure}
\begin{center}
\resizebox{1.15in}{1.15in}{\includegraphics{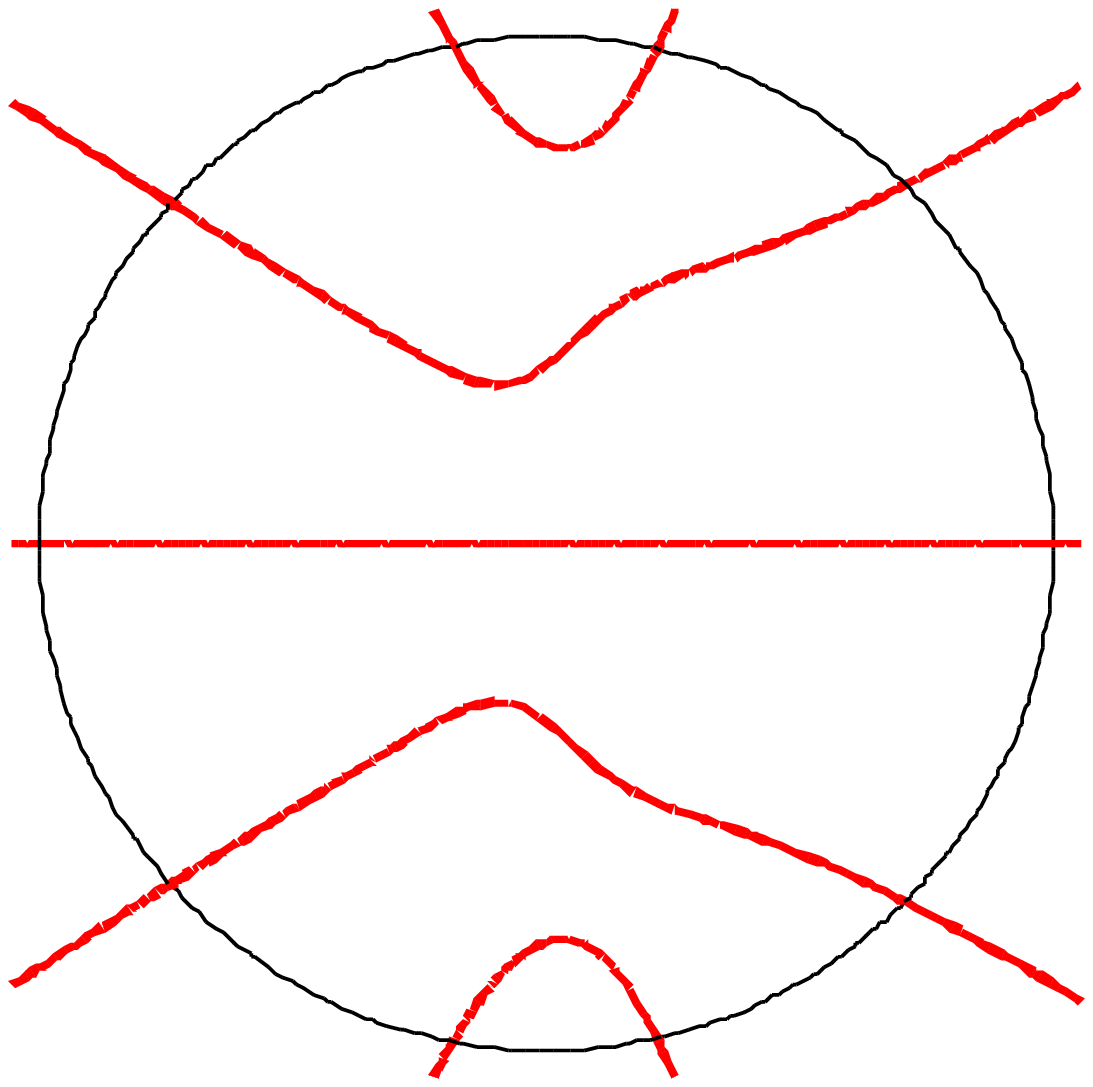}}
\hfill
\resizebox{1.15in}{1.15in}{\includegraphics{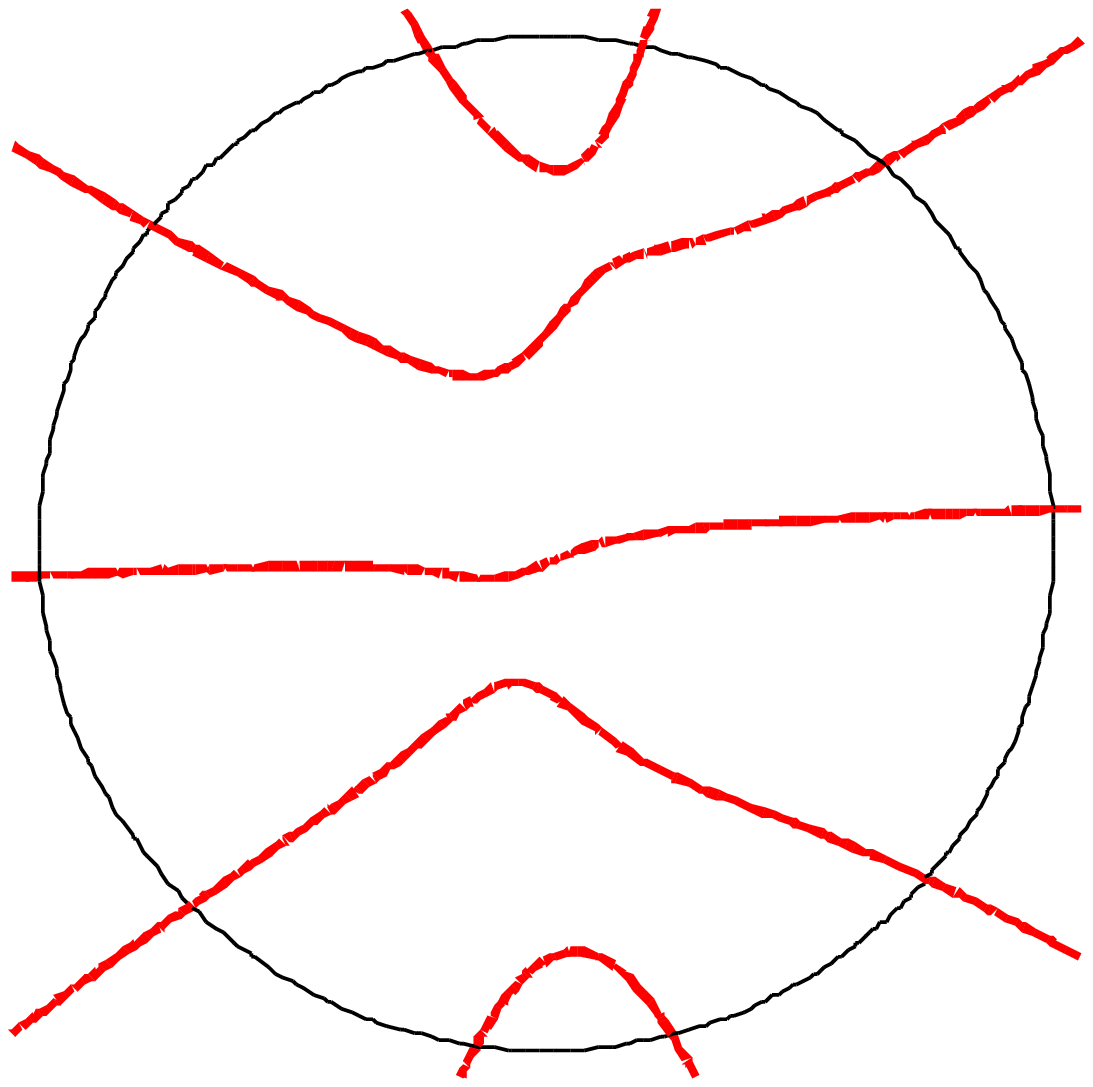}}
\hfill
\resizebox{1.15in}{1.15in}{\includegraphics{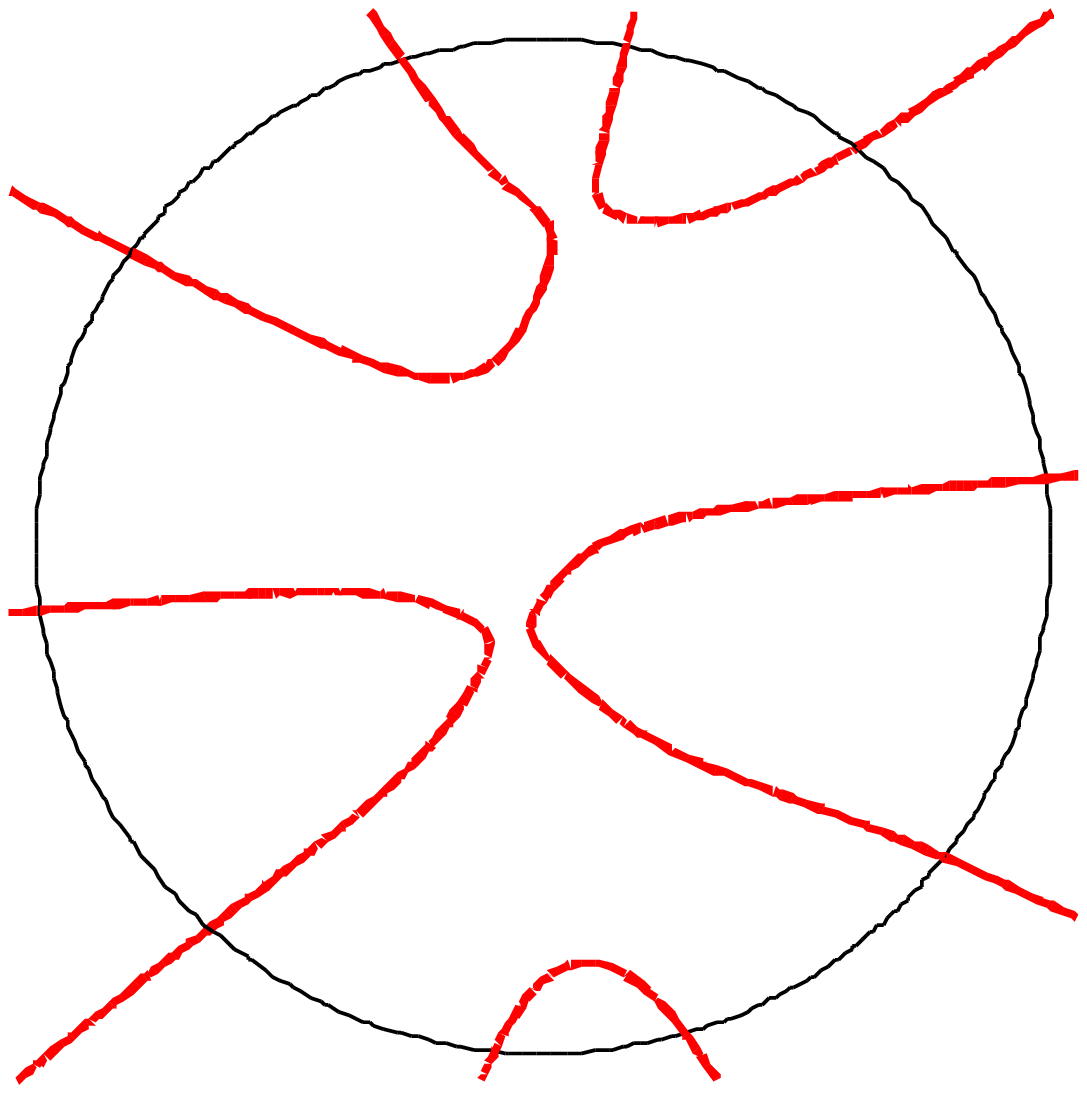}}
\hfill
\resizebox{1.15in}{1.15in}{\includegraphics{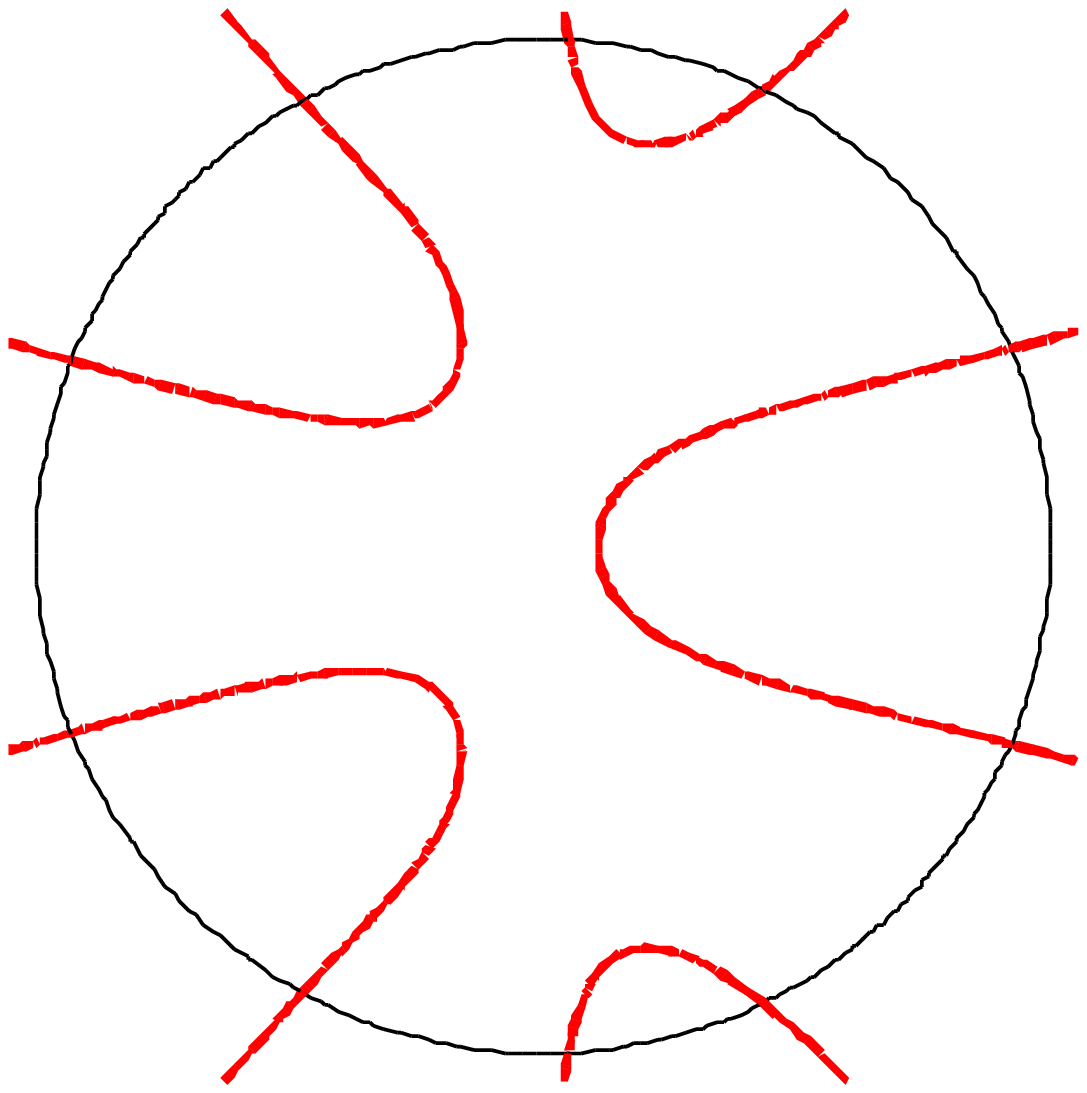}}
\caption{The curves $\Ct(f)$, where $f(z)$ is the quintic
of Figure~\ref{basketball-example-figure} and $\theta=0,
\frac{\pi}{12},\frac{\pi}{6},\frac{\pi}{2}$ respectively.
\label{Ctheta-figure}}
\end{center}
\end{figure}

\begin{lemma}
The curve $\Ct(f)$ is singular if and only if there exists $z \in \Ct(f)$ such that $f'(z)=0$.
\end{lemma}

\begin{proof}
If $f(z) = u(x,y) + i v(x,y)$ with $u,v$ real, then $\Ct(f)$ is the curve
defined by $v(x,y) \cos\theta - u(x,y) \sin\theta = 0$.  It follows that $(x_0,y_0) \in
\Ct(f)$ is a singular point if and only if $\frac{\partial v}{\partial x} \cos\theta -
\frac{\partial u}{\partial x} \sin\theta$ and $\frac{\partial v}{\partial y} \cos\theta -
\frac{\partial u}{\partial y} \sin\theta$ vanish at $(x_0,y_0)$. By the Cauchy-Riemann
equations, it follows that the partial derivatives of $u$ and $v$ with respect to $x$
and $y$ all vanish at $(x_0,y_0)$, and therefore $z_0 = x_0 + i y_0$ is a zero of $f'(z)$.
\end{proof}

\begin{cor} \label{cor:repeated}
If $f(z)$ has no repeated roots, then the curves $\Ct(f)$ are singular for at most $n-1$
values of $\theta$.  If $f(z)$ has a repeated root, then the curves $\Ct(f)$ are all singular.
\end{cor}

\begin{proof}
If $f'(z) = 0$ and $f(z) \neq 0$, then $z$ lies on $\Ct(f)$ for exactly one
value of $\theta \in \R/\pi\Z$, namely $\theta \equiv \arg(f(z)) \pmod{\pi}$.  On the other
hand, if $f'(z)=0$ and $f(z)=0$, then $z$ is a singular point on all $\Ct(f)$.
\end{proof}

We can now prove:

\begin{thm}[Fundamental Theorem of Algebra]
Every polynomial $f(z)$ of degree $n$ has exactly $n$ roots.
\end{thm}

\begin{proof} We may suppose that $f(z)$ is monic.
If $f(z)$ has a repeated root, then it certainly has a root;
by induction on the degree, $f(z)$ has exactly $n$ roots.

Assume, then, that $f(z)$ has no repeated roots.  Choose any two
different angles $\alpha,\beta \in \R/\pi\Z$ such that $C_\alpha(f)$ and
$C_\beta(f)$ are nonsingular; this is possible by Corollary
\ref{cor:repeated}.  Applying the argument which proved
Proposition \ref{prop:nonsing} to the curves $C_\alpha(f)$ and
$C_\beta(f)$, the theorem follows.
\end{proof}

The idea at the heart of Gauss's proof is that when $r$ is sufficiently
large, the $2n$ points of $C_\alpha(f) \cap S_r$ are paired via components of
$C_\alpha(f)$.  This associates to $f$ a purely combinatorial
object: a \emph{matching}, or partition into $n$ subsets of size two,
of the set $C_\alpha(f) \cap S_r$.
In a similar way, we can construct a matching on the set $C_\beta(f)
\cap S_r$, and the existence of $n$ distinct roots follows from the
combinatorial topology of these two matchings.  In the remainder of
this article, we study which (pairs of) such matchings can
arise from polynomials.  We begin by formalizing our study of the
matchings which arise from the curves $C_\theta(f)$.

Suppose that $C_\theta(f)$ is nonsingular.  By the
maximum principle for harmonic functions, the curve $C_\theta(f)$
cannot have any bounded connected components. It follows that
every connected component of $C_\theta(f)$ is diffeomorphic to
$\R$.  By B\'ezout's theorem, $C_\theta(f)$ cannot meet the
quadratic curve $S_r$ in more than $2n$ points. For $r$
sufficiently large, each connected component of $C_\theta(f)$ must
meet $S_r$ in at least two points, so we deduce that $C_\theta(f)$
has at most $n$ connected components.

Recall, however, that the curve $C_\theta(f)$ has $2n$ asymptotes,
with slopes at angles $(k\pi+\theta)/n$ for $k=0,\ldots,2n-1$.
Each connected component approaches at most two of these asymptotes,
and it follows that $C_\theta(f)$ has precisely $n$ connected
components, each diffeomorphic to $\R$, and each approaching two of
these asymptotes.  Moreover, if $r$ is sufficiently large then
each point of $C_\theta(f) \cap S_r$ lies on a different one of
these asymptotes, so that which pairs of points $C_\theta(f) \cap
S_r$ are joined by arcs in $C_\theta(f) \cap D_r$ is determined
entirely by which pairs of asymptotes lie on the same connected
component of $C_\theta(f)$.

If $C$ is any real plane curve, write $C(r) = C \cap D_r$.
Suppose that $C(r)$ is nonsingular and has $m$ connected
components, each with precisely two points lying on the circle
$S_r$.  Label these points counterclockwise from $0$ to $2m-1$,
starting from the positive real axis. We make the following
exception to the labeling rule: if $C=C_0(f)$ for some polynomial
$f$ and $r$ is sufficiently large that $C(r)$ meets the component
of $C$ asymptotic to the positive real axis, then we begin
labeling from the point of $C(r)$ furthest along this asymptote.

Then $C(r)$ induces a matching $\Match(C,r)$ on the set
$\{0,\ldots,2m-1\}$, where $a,b$ are matched if and only if the
points labeled $a$ and $b$ lie on the same connected component of
$C(r)$.  If $C = \Ct(f)$ is nonsingular and $r$ is sufficiently
large, then we have observed that the matching $\Match(\Ct(f),r)$
on $\{0,\ldots,2n-1\}$ does not depend on $r$; we denote this
matching by $\Match(f,\theta)$.  Our object is to study the
matchings $\Match(f,\theta)$, as well as the pairs of matchings
$(\Match(f,\alpha),\Match(f,\beta))$.

\begin{remark}
Recall that $f$ is assumed to be monic, so that if $C = \Ct(f)$
and $\theta \in [0,\pi)$, then for $r$ sufficiently large the
point of $C \cap S_r$ which is labeled $k$ lies on the asymptote
at angle $\frac{k\pi+\theta}{n}$.
\end{remark}

\begin{remark} \label{rmk:tangency}
If $C = \Ct(f)$ is nonsingular, we say a few words for future
reference about when $\Match(C,r)$ is defined, i.e., when
each connected component of $C(r)$ has precisely two points
lying on the circle $S_r$.  Each connected component of $C$
is diffeomorphic to $\R$; since a connected component of
$C(r)$ is a bounded connected subset of $C$, it is
either a single point, or diffeomorphic to the interval $[0,1]$.
Certainly each connected component which passes through the
interior of the disk $D_r$ must intersect the boundary at least
twice.  It follows that $\Match(C,r)$ is well-defined as long as
$C$ has no point of tangency (either external or internal) to
the circle $S_r$.
\end{remark}

\begin{remark} \label{rmk:discriminant}
When $f(z)=z^2+bz+c$ is quadratic, the pair of matchings
$(\Match(f,0), \Match(f,\pi/2))$ is determined by the quadrant of the complex
plane containing the discriminant $b^2-4c$.
An analogous description must exist for polynomials of higher
degree.  Indeed, regarding a monic polynomial of degree $n$ as a vector in
$\R^{2n}$ by taking the real and imaginary parts of its
coefficients, the subset of $\R^{2n}$ for which at least one of the
curves $R$ and $I$ is singular is an algebraic set. Its complement,
the set of polynomials for which $\Match(f,0)$ and $\Match(f,\pi/2)$ are both
defined, is therefore a real semialgebraic set.  However, every
connected component of a real semialgebraic set is again
semialgebraic \cite[Thm.~2.4.5]{Bochnak}; since the set of
polynomials which yield a particular pair of matchings is a union of
such connected components, that set is itself semialgebraic.
\end{remark}

%===========================================================================

\section{Combinatorics of basketballs} \label{section:basketball}

Throughout, if $a\leq b$ are integers,
we write $[a,b]$ for the set $\{a,a+1,\dots,b-1,b\}$.  (It should be clear from context
whether this notation refers to an interval in $\Z$ or in $\R$.)

\begin{defn} Let $a\leq b$ be integers.
A \defterm{partition} of $[a,b]$ is a collection of pairwise disjoint sets (called
\defterm{blocks}) whose union is $[a,b]$. The number of blocks is the \defterm{order} of
the partition. Two blocks are said to \defterm{cross} if there are integers
$i<j<k<\ell$ such that $i,k$ belong to one block and $j,\ell$ belong to the other block.
If no two blocks cross, then the partition is said to be \defterm{noncrossing}. A
\defterm{matching} is a partition in which every block has cardinality 2.
\end{defn}

\begin{figure}[h]
\begin{center}
\resizebox{2.5in}{1in}{\includegraphics{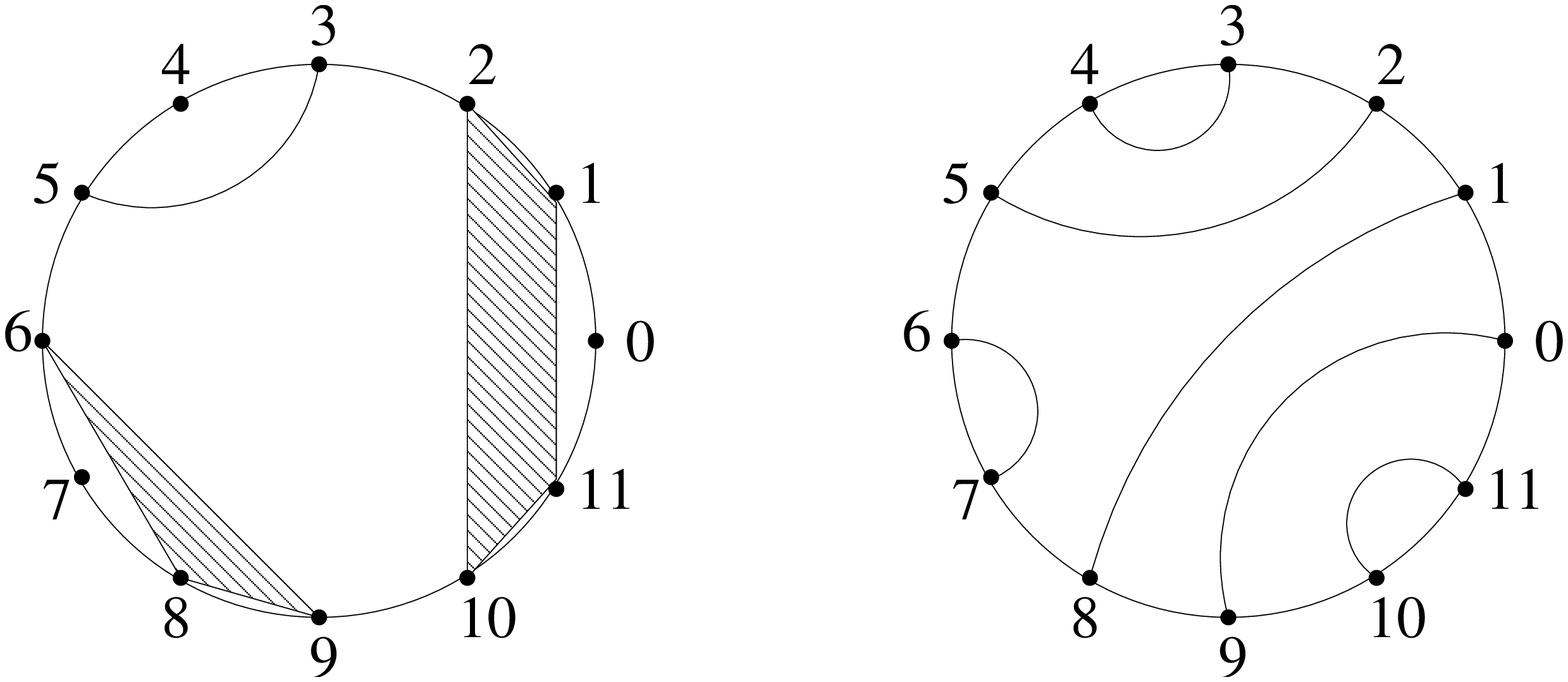}}
\caption{\emph{Left}: The noncrossing partition whose non-singleton blocks are
$\{1,2,10,11\}$, $\{3,5\}$ and $\{6,8,9\}$.  \emph{Right}: A noncrossing matching.
\label{NC-example-figure}}
\end{center}
\end{figure}

A noncrossing partition can be represented by
placing the numbers $a,a+1,\dots,b$ around a circle and connecting numbers in the
same block, as in Figure~\ref{NC-example-figure}.

\begin{example}  If the curve $C_\theta(f)$ is nonsingular, then
the matching $\Match(f,\theta)$ is
noncrossing by the Jordan Curve Theorem.
\end{example}

The theory of noncrossing partitions comprises a substantial chapter in modern combinatorics;
see \cite{Simion} for a comprehensive survey.  In the present study, we will be concerned
most with noncrossing matchings on the sets
  $$\En=\{0,2,4,\dots,4n-2\}\qquad\text{and}\qquad \On = \{1,3,\dots,4n-1\},$$
where $n$ is a positive integer.

It is well-known (and elementary to verify) that the number of matchings on $\En$
(resp.\ $\On$) is $(2n-1)(2n-3)\cdots(3)(1) = \frac{(2n)!}{2^nn!}$, and that the number of noncrossing
matchings is the Catalan number $\frac{1}{n+1}\binom{2n}{n}$.

\begin{lemma} \label{lemma:parity-NCM}
Let $M$ be a matching on $[0,2n-1]$.  If $M$ is noncrossing, then exactly one member of each pair
$\{i,j\}$ is even.
\end{lemma}

\begin{proof}
Suppose that $i < j$ and $i\equiv j \pmod{2}$.  Then the sets
  $$X = [i+1,j-1], \qquad Y = [0,i-1] \cup [j+1,2n-1]$$
both have odd cardinality, so some $x\in X$ is paired with some $y\in Y$. Then $\{x,y\}$
crosses $\{i,j\}$, contradicting the condition that $M$ is noncrossing.
\end{proof}

Equivalently, if $M$ is a noncrossing matching on $\En$ or on $\On$, then every pair
in $M$ consists of two numbers that are non-congruent modulo~4.

\begin{defn}
A \defterm{bimatching of order $n$} is a pair $B=(\Beven,\Bodd)$,
where $\Beven$ is a matching on $\En$ and $\Bodd$ is a matching on
$\On$. A pair in $\Beven$ (resp.\ $\Bodd$) is called an
\defterm{even} (resp.\ \defterm{odd}) pair of $B$.  Let
$\Meven,\Modd$ denote the matchings induced on $[0,2n-1]$ by
$\Beven,\Bodd$ via the maps $[0,2n-1] \rightarrow \En,\On$ sending $i
\mapsto 2n,2n+1$ respectively.  We say that $B$ is the bimatching
corresponding to the ordered pair of matchings $(\Meven,\Modd)$,
and vice-versa.
\end{defn}

%It is immediate that the bimatchings $B=(\Beven,\Bodd)$
%of order $n$ in which both $\Beven$ and $\Bodd$ are noncrossing
%are counted by the square of the $n$th Catalan number.

\begin{example} \label{ex:bimatch} Choose $\alpha, \beta \in [0,\pi)$ such that
$C_{\alpha}(f)$ and $C_{\beta}(f)$ are nonsingular, and suppose
that $\alpha < \beta$.  Then we obtain a bimatching $\Basket(f,\alpha,\beta)$
on $[0,4n-1]$ corresponding to the pair of matchings
$(\Match(f,\alpha),\Match(f,\beta))$,

This bimatching can be obtained geometrically as follows.  For~$r$
sufficiently large, label the points of $(C_\alpha(f) \cup
C_\beta(f)) \cap S_r$ counterclockwise from $0$ to $4n-1$,
beginning with the point of $C_\alpha(f)$ lying on the asymptote
at angle $\alpha/n$. Again for~$r$ sufficiently large, each point
of $C_\alpha(f) \cap S_r$ will be labeled with an element of
$\En$, and each point of $C_\beta(f) \cap S_r$ will be labeled with
an element of $\On$.  The bimatching is then induced by the arcs
of $C_\alpha(f)\cap D_r$ and $C_\beta(f)\cap D_r$ respectively.  See, e.g.,
Figure~\ref{basketball-example-figure}.
\end{example}

\begin{lemma} \label{lemma:odd-cross}
Let $B=(\Beven,\Bodd)$ be a bimatching of order $n$. Suppose that
both $\Beven$ and $\Bodd$ are noncrossing.  Then each pair of
$\Beven$ (resp.\ $\Bodd$) crosses an odd number of pairs of
$\Bodd$ (resp.\ $\Beven$).
\end{lemma}

\begin{proof}
Let $i,j$ be an odd pair of $B$ with $i<j$.  Let
  $$X=\{i+1,i+3,\dots,j-1\}, \qquad
    Y=\{0,2,\dots,i-1,\; j+1,\dots,4n-2\}.$$
Then a pair $k,\ell\in \Beven$ crosses $i,j$ if and only if exactly one of $k,\ell$ belongs
to $X$. By Lemma~\ref{lemma:parity-NCM}, we have $i\not\equiv j\pmod{4}$, which implies that
$\#X$ and $\#Y$ are both odd. Therefore there are an odd number of such pairs $k,\ell$. By the
same argument, each pair in $\Beven$ crosses an odd number of pairs in $\Bodd$.
\end{proof}

In particular, each pair of $B$ crosses at least one pair of the opposite parity.  We are
interested primarily in the case when no extra crossings occur.

\begin{defn}
A \defterm{basketball of order $n$} (or simply an
\defterm{$n$-basketball}) is a bimatching $B=(\Beven,\Bodd)$ of
order~$n$ in which the matchings $\Beven,\Bodd$ are both noncrossing,
and each pair $e\in\Beven$ crosses exactly one pair $o\in\Bodd$.
The pair of pairs $e,o$ is called a \defterm{quartet}. The set
of all $n$-basketballs is denoted $\BB_n$.
\end{defn}

\begin{figure}[h]
\begin{center}
\hfill\phantom{.}
\resizebox{1.5in}{1.5in}{\includegraphics{basketball1.eps}}
\hfill
\resizebox{1.5in}{1.5in}{\includegraphics{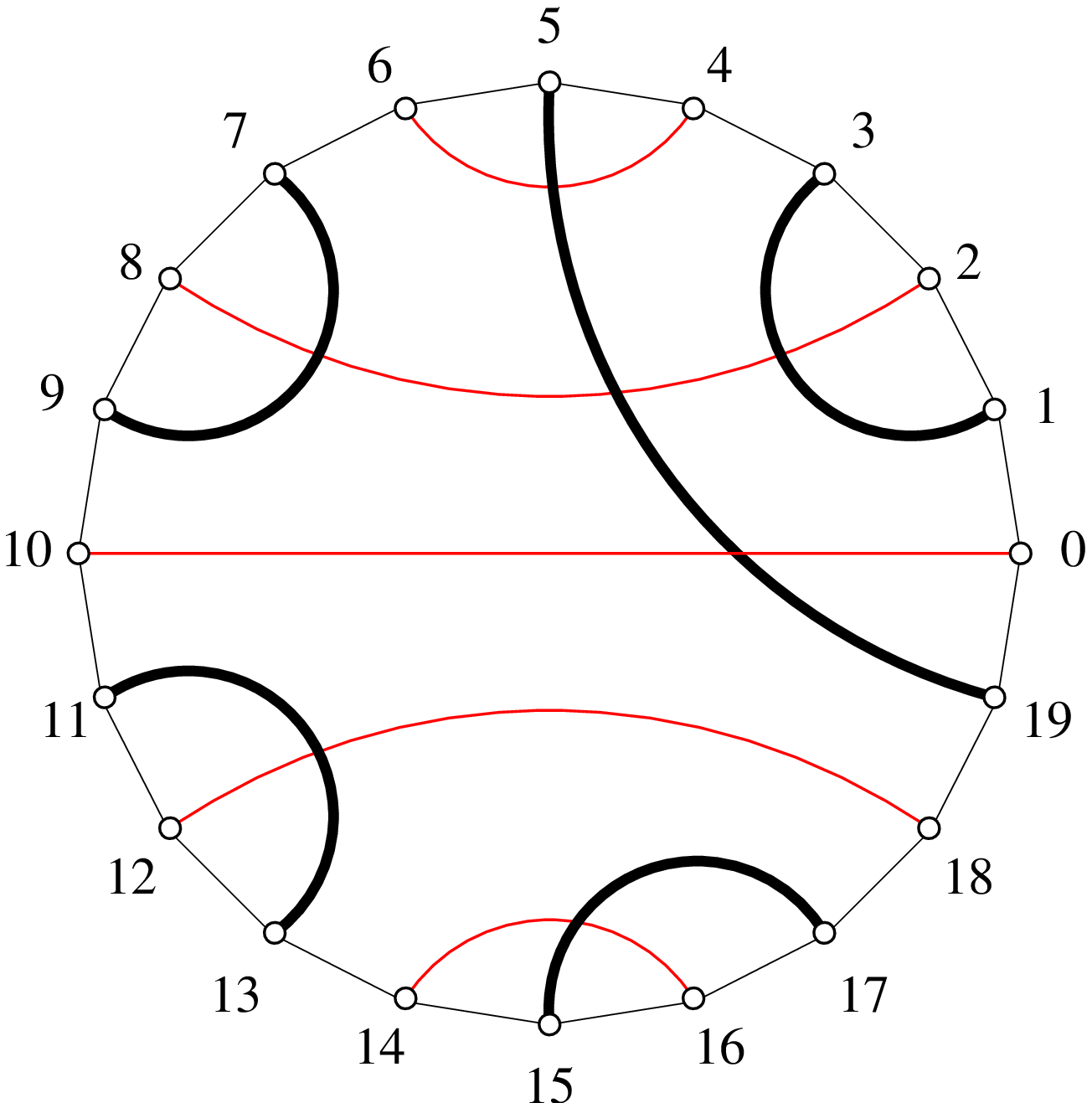}}
\hfill\phantom{.}
\end{center}
\caption{
  \emph{Left}: A basketball of order 5.
  \emph{Right}: A bimatching that is not a basketball
    (the pair $\{5,19\}$ crosses three even pairs).
\label{bimatching-example-figure}}
\end{figure}

The bimatching $\Basket(f,\alpha,\beta)$ of Example~\ref{ex:bimatch} is a
basketball of order $n$.  Indeed, every crossing between an even
pair (i.e., an arc of $C_\alpha(f)$) and an odd pair (i.e., an arc
of $C_\beta(f)$) corresponds to a root of $f$, and there are
exactly $n$ roots.

\begin{thm} \label{thm:count-bballs}
The number of basketballs of order $n$ is $\frac{1}{3n+1}\binom{4n}{n}$.
\end{thm}

\begin{proof}
By a special case of \cite[Lemma~4.1]{Edelman}, the given quantity enumerates the set $\QQ_n$
of noncrossing partitions of $[0,4n-1]$ into $n$ blocks of size 4. We will exhibit a bijection
$\phi:\BB_n\to\QQ_n$.

Let $B \in \BB_n$.  Define a partition $\phi(B)$ of $[0,4n-1]$ into $n$ blocks of size 4 by
replacing each quartet $\{i,j\},\{k,\ell\}$ with the block $\{i,j,k,\ell\}$. It is elementary
to deduce from the definition of a basketball that $\phi(B)$ is noncrossing, hence belongs
to $\QQ_n$.  It is also evident that the function $\phi$ is injective.

Now let $Q\in\QQ_n$ and let $K=\{a<b<c<d\}$ be a block of $Q$.  Since $Q$ is noncrossing, each
of the sets
  \begin{equation} \label{inbetween}
    [a+1,b-1],\quad
    [b+1,c-1],\quad
    [c+1,d-1],\quad
    [d+1,4n-1]\cup[0,a-1]
  \end{equation}
must be a union of blocks, hence must have cardinality divisible
by 4. In particular, $a\equiv c\not\equiv b\equiv d \pmod{2}$. It
follows that replacing $K$ with the two pairs $\{a,c\},\{b,d\}$,
and doing the same for every other block of $Q$, yields a
$n$-basketball $B$ such that $\phi(B)=K$.  Therefore $\phi$ is
surjective.
\end{proof}

\begin{figure}[h]
\begin{center}
\resizebox{3.6in}{1.5in}{\includegraphics{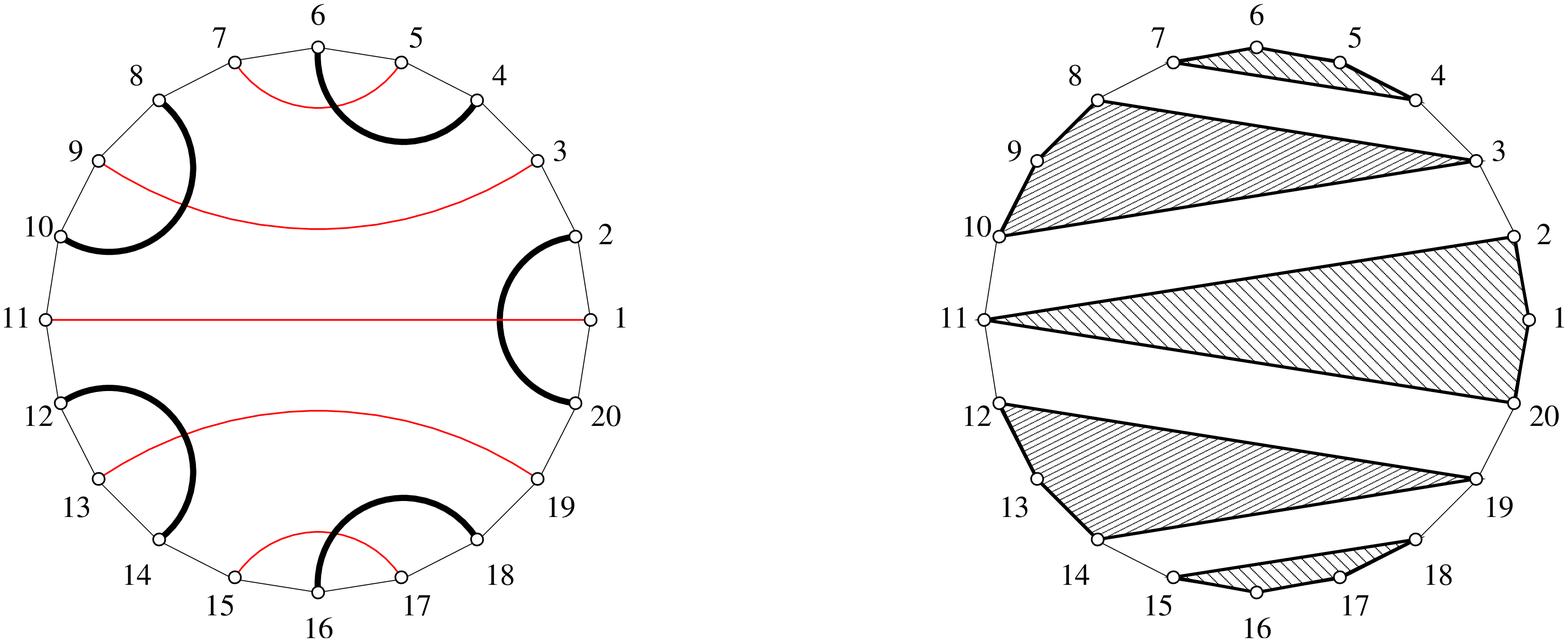}}
\end{center}
\caption{The bijection $\phi$ of Theorem~\ref{thm:count-bballs},
mapping basketballs of order $n$ (\emph{left}; here $n=5$)
to noncrossing partitions of $[0,4n-1]$ with all blocks of size 4 (\emph{right}).}
\label{bijection-figure}
\end{figure}

The numbers $\frac{1}{3n+1}\binom{4n}{n}$ form sequence A002293 in the
On-Line Encyclopedia of Integer Sequences \cite{EIS}.
Other combinatorial objects enumerated in the same way include
  \begin{itemize}
  \item quaternary trees with~$n$ internal nodes;
  \item dissections of a $(3n+2)$-gon into $n$~pentagons; and
  \item rooted plane maps (see \cite{LW}).
  \end{itemize}

\begin{defn}
Let $B$ be a basketball.  An \defterm{ear} of $B$ is a quartet $\{i,j\},\{k,\ell\}$ such that
the numbers $i,j,k,\ell$ are consecutive (in some order) modulo $4n$.
\end{defn}

For instance, the basketball shown in Figure~\ref{bijection-figure} has two ears: $\{4,6\},\{5,7\}$
and $\{15,17\},\{16,18\}$.

\begin{prop} \label{thm:ear-thm}
Every basketball $B$ of order $n\geq 2$ has at least two ears.
\end{prop}

\begin{proof}
We proceed by induction on $n$.  The base case $n=2$ is amenable to proof by inspection, as
there are only four basketballs of order two.

Suppose now that $n>2$.  Let $B$ be a basketball and
let $T=\{i,j\},\{k,\ell\}$ be a quartet with $i<j<k<\ell$..  If $T$ is an ear, then $B-T$ is
a basketball on the ground set
$[\ell+1,4n-1]\cup[0,i-1]$. By induction, $B-T$ contains
at least two ears $E,F$.  Without loss of generality, we have
$\ell+1\not\in E$.  Then $E$ is an ear of $B$.

Now suppose that $T$ is not an ear.  Then at least two of the four sets of \eqref{inbetween}
are nonempty.  By a similar argument, each of those two sets contains an ear that is also an
ear of $B$.
\end{proof}

We point out several combinatorial problems related to bimatchings and basketballs.
\begin{enumerate}
\item
The set of $n$-basketballs is invariant under the operation
of \emph{rotation}: replacing each quartet $\{i,j\},\{k,\ell\}$ with
$\{i-2,j-2\},\{k-2,\ell-2\}$, with all indices taken modulo $4n$.
Let $r(n)$ be the number of $n$-basketballs up to rotation, so that
$(r(1),r(2),\dots)=(1,2,6,22,103,614,3872,\dots)$.  These numbers occur as sequence
A103941 in \cite{EIS}, enumerating
unrooted loopless plane maps with $n$ edges \cite[Theorem~4.4]{LW}.
There does not seem to be an obvious relationship between plane maps
and basketballs.
\item  One might instead seek to count equivalence classes of
basketballs up to \emph{half-rotation}, replacing $\{i,j\},\{k,\ell\}$
with $\{k-1,\ell-1\},\{i-1,j-1\}$, or up to half-rotation and \emph{reflection},
replacing $\{i,j\},\{k,\ell\}$ with $\{-i,-j\},\{-k,-\ell\}$.
\item It would be interesting to enumerate the
bimatchings of order $n$ by total number of crossings.
\end{enumerate}

%===========================================================================

\section{The Inverse Basketball Theorem} \label{section:IBT}

In light of Gauss's proof of the FTA, it is natural to ask whether
every basketball of order~$n$ arises as $B(f,0,\Pi/2)$ for some
suitably chosen polynomial $f(z)$ of degree $n$.  The main result of
this section is that something more general is true.

\begin{thm} \label{IBT}
Fix any $\alpha,\beta \in [0,\pi)$ with $\alpha < \beta$.
Then every basketball has the form $\Basket(f,\alpha,\beta)$
for some monic polynomial $f$.
\end{thm}

Before plunging into the details, let us describe the argument
informally. Proposition~\ref{thm:ear-thm} suggests an inductive
approach. Given an $n$-basketball $B$ for which we want to
construct a realizing polynomial, we would like to remove an ear
from $B$, inductively construct a realizing polynomial $f(z)$ for
the resulting $(n-1)$-basketball $B'$, and then insert the missing
ear by replacing $f(z)$ with $f(z)(z-R)$ for some suitable
$R$. If we choose $R$ to be much greater in absolute value than
any of the roots of $f$, then the components of $B'$ will be
perturbed only slightly, and will retain their combinatorial
structure near the origin, where $\arg(z-R)\approx \pi$ is
close to $0 \pmod{\pi}$.  (No confusion
should arise between the root $R$ and the curve $R=\{z\st\Re f(z)=0\}$,
which does not appear in this section.)

\begin{example}
\label{cubic-example}
 Let $f(z)=z^3+iz^2+z-2$, a cubic whose associated
basketball is shown on the left of Figure~\ref{insert-example-figure}.
Every root of $f(z)$ has complex magnitude $<2$.  If we choose a new root
that is much larger in magnitude, say $R=8+8i$, then the basketball
of $f(z)(z-R)$ is given combinatorially by ``inserting an ear
at $R$'', as shown.

\begin{figure}[h]
\begin{center}
\hfill\phantom{.}
\resizebox{1.5in}{1.5in}{\includegraphics{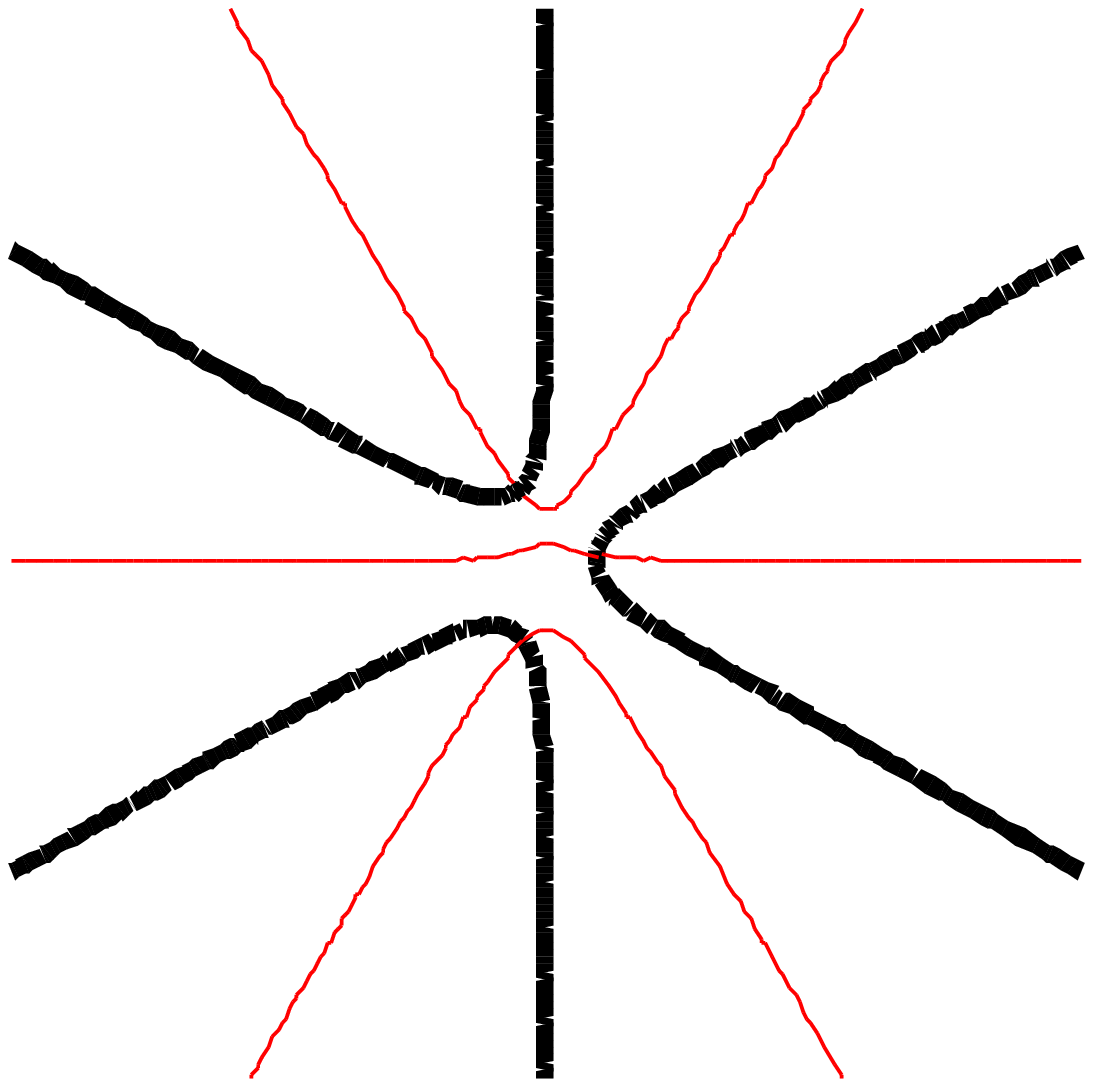}}
\hfill
\resizebox{1.5in}{1.5in}{\includegraphics{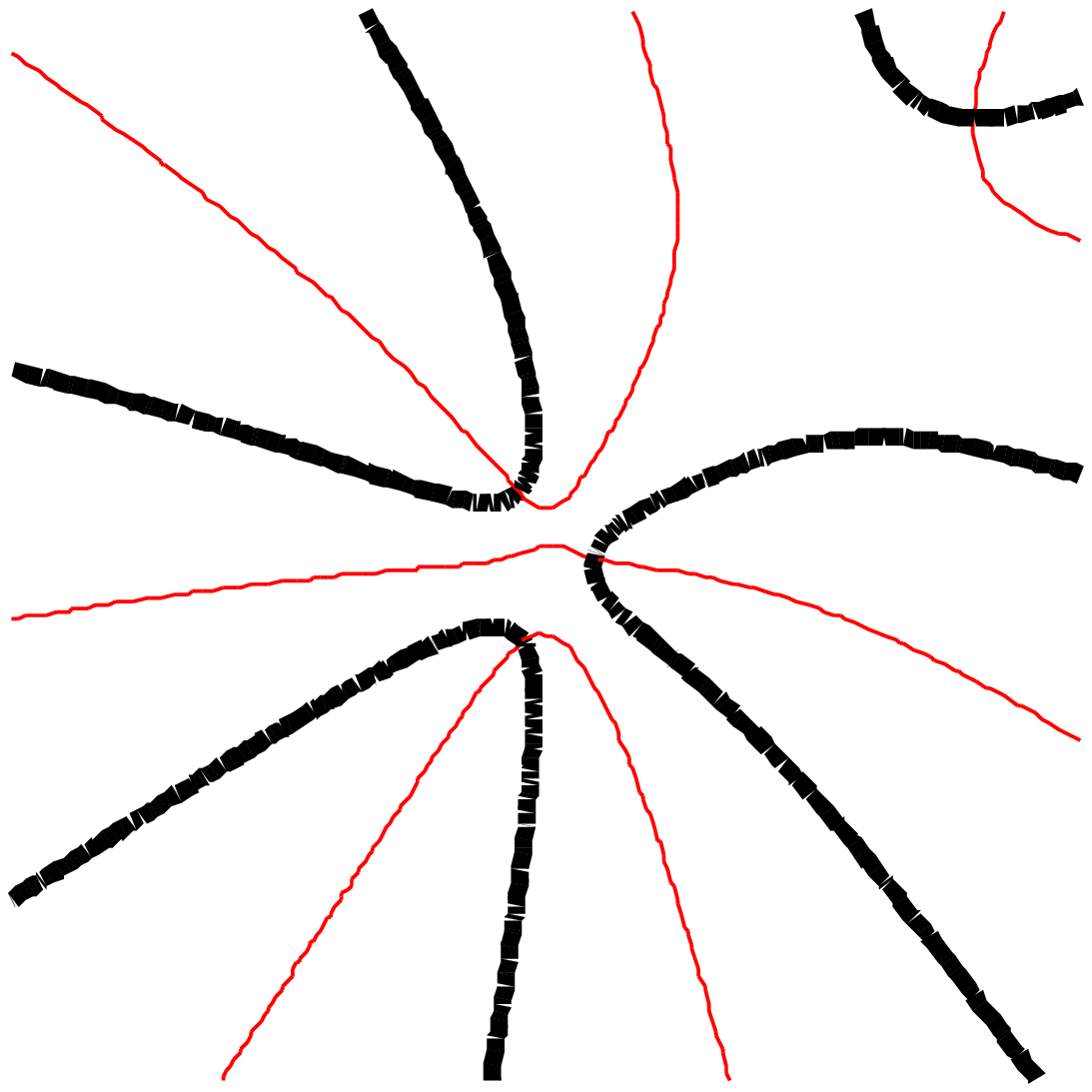}}
\hfill\phantom{.} %% this is a kludge to get LaTeX to space things out nicely
\end{center}
\caption{
  The basketballs of $f(z)$ (\emph{left}) and of $f(z)(z-R)$ (\emph{right}),
  where $f(z)$ is the cubic of Example~\ref{cubic-example}.
\label{insert-example-figure}}
\end{figure}
\end{example}

The following theorem explains how
$\Match(f(z)(z-R),\theta)$ is obtained from $\Match(f,\theta)$, and is the
crucial result from which Theorem~\ref{IBT} will follow.  We first introduce some
useful notation.  If $M$ is a matching on $[0,2n-1]$, define a matching $\widehat{M}$ on $[0,2n+1]$ by
  $$\widehat{M} = \big\{ \{0,2n+1\} \big\} \;\cup\; \big\{ \{i+1,j+1\} \st \{i,j\} \in M \big\}.$$
Also, for any positive real $R$, let $g_R(z) = (z-R)f(z)$.

\begin{thm} \label{thm:matchings}
Let $f(z)$ be a monic polynomial, and let $\theta \in (0,\pi)$ such that $\Ct(f)$ is
nonsingular.  Then, for $R$ sufficiently large, we have
  $$\Match(g_R,\theta) = \widehat{\Match(f,\theta)}.$$
\end{thm}

Before we can prove Theorem~\ref{thm:matchings}, we must develop
several subsidiary results.  The first of these is the following,
which is implicit in the statement of Theorem~\ref{thm:matchings}.

\begin{prop} \label{prop:smooth}
Under the assumptions of Theorem~\ref{thm:matchings}, the curves
$\Ct(g_R)$ are nonsingular for $R$ sufficiently large.
\end{prop}

To establish this fact, we first need the following technical result.
The idea is that as $R\to\infty$, one root of
$g_R'(z)$ increases without bound, while the other $n-1$ roots approach
the roots of $f'(z)$.

\begin{lemma} \label{lemma:singularities}
Let $\alpha_1,\ldots,\alpha_{n-1}$ denote the roots of $f'(z)$.  For $R\in\R$ sufficiently large,
there is an ordering $\beta_1(R),\ldots,\beta_n(R)$ of the roots
of $g_R'(z)$ with the following properties:\footnote{
  For simplicity of notation, we will often abbreviate $\beta_i(R)$ by $\beta_i$.}

  \begin{enumerate}
    \item if $i \le n-1$, then $\beta_i \rightarrow \alpha_i$ as $R \rightarrow \infty$;
    \item $|\beta_n| \rightarrow \infty$ and $\arg(\beta_n)\to 0$ as $R \rightarrow \infty$;
    \item if $i \le n-1$, then $\arg(g_R(\beta_i)) \rightarrow \arg(f(\alpha_i)) + \pi$ as
        $R \rightarrow \infty$;
    \item $\arg(g_R(\beta_n)) \rightarrow \pi$ as $R \rightarrow \infty$.
  \end{enumerate}
\end{lemma}

\begin{proof}
Note that $g_R'(z) = f(z) + (z-R)f'(z)$.  Put $t = 1/R$, so that
  $$g_R'(z)/R = (tz-1)f'(z) + tf(z).$$
Note that $\lim_{t\to 0} g'_R(z)/R=-f'(z)$.  The existence of an
ordering of the roots with property (1) now follows from the argument
principle~\cite[Theorem~18,~p.~152]{Ahlfors}.
Assertion (3) is evident: we have seen that~$\beta_i(R)$ is
bounded, and therefore that $\arg(\beta_i(R)-R) \rightarrow \pi$ as $R\to\infty$.

For (2) and (4), observe from the
coefficients of $z^n$ and $z^{n-1}$ in $g_R'$ that the sum of the roots of $g_R'(z)$ is
  $$\frac{n}{n+1} (R - a_{n-1}),$$
where $a_{n-1}$ is the coefficient of $z^{n-1}$ in $f(z)$.
Since the sum $\beta_1+\cdots+\beta_{n-1}$ is bounded, so is
  $$\beta_n(R) - \frac{n}{n+1} R = -\left(\frac{n}{n+1}a_{n-1}+\beta_1+\cdots+\beta_{n-1}\right).$$
That is, $\beta_n(R)$ is within a bounded distance of $\frac{n}{n+1}R$, and so we have
  $$\arg(\beta_n) \rightarrow 0, \qquad \arg(\beta_n-R) \rightarrow \pi, \qquad \text{and} \qquad
   |\beta_n|\rightarrow \infty,$$
whence
  $$\arg(g_R(\beta_n)) = \arg((\beta_n-R)f(\beta_n)) \rightarrow \pi \,.$$
\end{proof}

\begin{proof}[Proof of Proposition~\ref{prop:smooth}]
For $R$ sufficiently large, $g_R(z)$ has no repeated roots.  By
Lemma \ref{lemma:singularities}, we see that
$\arg(g_R(\beta_i(R)))$ becomes arbitrarily close to
$\arg(f(\alpha_i))+ \pi$ for $i<n$, and to~$\pi$ for $i=n$.  Using
the hypotheses that $\Ct(f)$ is nonsingular and $\theta \neq 0$,
it follows that $\arg(g_R(\beta_i(R)))$ is not congruent to
$\theta \pmod{\pi}$ for any $i$, and so $C_{\theta}(g_R)$ is
nonsingular for $R$ sufficiently large.
\end{proof}

The next step in proving Theorem~\ref{thm:matchings} is to show that
the topological behavior of the curve $\Ct(g_R)$ inside the disk $D_r$
is the same as that of the curve $\Ct(f)$; that is, the passage from
$f(z)$ to $g_R(z)$ preserves the combinatorial
type of the matching $\Match(f,\theta)$.

\begin{prop} \label{prop:disk}
Let $f(z)$ be a monic polynomial of degree $n$, let $\theta \in
(0,\pi)$ be such that $\Ct(f)$ is nonsingular, and choose $r$
so that $\Ct(f)$ is not tangent to the circle $S_r$.  (Hence
$\Match(\Ct(f),r)$ is well-defined by Remark~\ref{rmk:tangency}.)
Then for $R$ sufficiently large, we have
  $$\Match(\Ct(g_R),r) = \Match(\Ct(f),r) \,.$$
Moreover, as $R\rightarrow \infty$, the
point labeled $i$ on $\Ct(g_R) \cap S_r$ approaches the point
labeled~$i$ on $\Ct(f) \cap S_r$.
In particular, if $r$ is sufficiently large, then for $R$ sufficiently large we have
  $$\Match(\Ct(g_R),r) = \Match(f,\theta)\,.$$
\end{prop}

In order to prove Proposition~\ref{prop:disk},
we will need the following (presumably standard) facts from metric topology.

\begin{lemma} \label{lemma:technical}
Let $(X,d)$ be a compact, locally connected metric space.  For $x \in
X$ and $\epsilon > 0$, write $B(x,\epsilon)$ for the open ball of
radius $\epsilon$ centered at $x$.  For any subset $K$ of $X$ and
any $\epsilon > 0$, define
  $$N_{\epsilon}(K) = \bigcup_{x \in K} B(x,\epsilon) \,.$$
For $\beta>0$ and $t \in (0,\beta)$, let $h_t : X\rightarrow \R$
be a family of continuous functions converging pointwise (hence uniformly)
to some continuous $f : X \rightarrow \R$ as $t \rightarrow 0$.  Define
  \begin{align*}
  Z &= \{x\in X \st f(x)=0\},\\
  Y_t &= \{x\in X \st h_t(x)=0\}.
  \end{align*}
Then:
  \begin{enumerate}
    \item For all $\epsilon > 0$, if $t$ is sufficiently close to $0$ then $Y_t \subset
        N_\epsilon(Z)$.
    \item Suppose furthermore that $X\subset\R^2$, that $U\subset X$ is an open
        set on which the functions $f, h_t$ are harmonic, and that $f$ is nowhere locally identically
        zero on $U$. Then, for all $\epsilon > 0$, there exists $t_0=t_0(\epsilon)>0$
        such that $Z \cap U \subset N_{\epsilon}(Y_t \cap U)$ for all $t<t_0$.
  \end{enumerate}
\end{lemma}

\begin{proof}
Suppose that (1) fails: that is, there exist $\epsilon > 0$ and sequences $t_n
\rightarrow 0$, $y_n \in Y_{t_n}$ such that for all $n$ we have $B(y_n, \epsilon) \cap Z =
\emptyset$.  Since $X$ is compact, the sequence~$\{y_n\}$ has a limit point; therefore,
replacing~$y_n$ and~$t_n$ by suitable subsequences, we may assume without loss of generality
that the sequence $y_n$ converges to a point $y$.

Note that $y \not\in Z$.  Indeed, $B(y,\epsilon/2) \subset B(y_n,\epsilon)$ for $n$
sufficiently large, so in fact $B(y,\epsilon/2) \cap Z = \emptyset$.  Choose a connected
neighborhood $U \subset B(y,\epsilon/2)$ of $y$; by the Intermediate Value Theorem, $f$ is
either strictly positive or strictly negative on $U$.  Since $f$ is continuous and $X$ is
compact, we find a neighborhood $V \subset U$ of $y$ and a constant $\delta > 0$ such that
  \begin{equation} \label{eq:delta}
    |f(x)| > \delta > 0 \qquad \text{for all} \ x \in V\,.
  \end{equation}
Now, choose $n$ sufficiently large so that $y_n \in V$ and (by uniform convergence) so that
$|g_{t_n}(x) - f(x)| < \delta/2$ for all $x \in X$.  Then $f(y_n) > \delta$ by
\eqref{eq:delta}.  On the other hand, $g_{t_n}(y_n) = 0$, so $|f(y_n)| = |g_{t_n}(y_n) -
f(y_n)| < \delta/2$.  This is a contradiction.

To prove (2), fix $\epsilon > 0$ and cover $Z \cap U$ with finitely many open balls
$B_1,\ldots,B_n$ of radius $\epsilon/2$.  (To do so, first cover $Z$ by
finitely many such balls, and then discard the ones that do not meet $Z \cap U$.)

For each $i = 1,\ldots,n$, the ball $B_i$ contains a point $z_i$ of $Z \cap U$.  By the
maximum principle \cite[Theorem~21, p.~166]{Ahlfors} applied in a neighborhood of $z_i$
contained in $B_i \cap U$, we can find points $p_i,q_i$ in $B_i \cap U$ such that $f(p_i) <
f(z_i) = 0 < f(q_i)$.  Choose $\delta > 0$ such that $\delta < |f(p_i)|,|f(q_i)|$ for all
$i$.

By uniform convergence, we can choose $t_0=t_0(\epsilon) > 0$ so that
  $$|h_t(x) - f(x)| < \frac{\delta}{2} \qquad \text{for all} \ t \in (0,t_0), \ x \in X\,.$$
In particular,
  $$h_t(p_i) < -\frac{\delta}{2} < 0 < \frac{\delta}{2} < h_t(q_i)$$
for all $i=1,\ldots,n$ and $t \in (0,t_0)$.  Since the function $h_t$ is continuous, it has a
root $w_i \in B_i \cap U$ for every $i$ (for example, along any path joining $p_i$ and $q_i$).

Now, for any $z \in Z \cap U$, we can find $i \in \{1,\ldots,n\}$ such that $z \in B_i$.
Since also $w_i\in B_i$, we have $|z-w_i|<\epsilon$ as desired.
\end{proof}

\begin{remark} Under the hypotheses of part (2) of
Lemma~\ref{lemma:technical}, the conclusions of part (1) and
(2) together imply that the \textit{Hausdorff distance} between
$Z$ and $Y_t$, namely $\max\left(
\inf(\epsilon>0\st Z\subset N_\epsilon(Y)),\inf(\epsilon>0\st Y\subset N_\epsilon(Z))
\right)$, tends to zero as $t\to 0$.  Of course, part (2) of Lemma~\ref{lemma:technical} holds for any
family of functions $\{h_t\},f$ that satisfy the maximum principle, even if they are not harmonic.
\end{remark}

\begin{proof}[Proof of Proposition \ref{prop:disk}]
By Lemma \ref{lemma:singularities}, the curves $\Ct(g_R)$ are nonsingular for $R$ sufficiently
large, say $R>1/\beta$.  Let $t = 1/R$, and write $H_t(z) = g_R(z)/R = (tz-1)f(z)$.  (Note that $H_t(z)$ is
well-defined for~$t=0$.)  Then $\{\Ct(H_t)\}_{t\in[0,\beta)}$ is a
family of nonsingular real algebraic curves for $t \in [0,\beta)$, and
$\Ct(H_t) = \Ct(g_R)$. Let~$m$ be the number of connected components
of $\Ct(f)$.  For $0 \le j
\le 2m-1$, let~$P_j$ denote the point of $\Ct(f)(r)$ that is labeled $j$.
Let $h_t(x,y)$ be the real polynomial $\Im(e^{-i\theta} H_t(x+iy))$.

We will eventually apply Lemma~\ref{lemma:technical} to the family~$\{h_t\}$
on the domain $X = D_r$.  Accordingly, we define
  \begin{align*}
    Y_t &\,=\, \Ct(H_t) \cap D_r \,=\, \Ct(g_R) \cap D_r,\\  %%    \Ct(H_t)(r)
    Z   &\,=\, \Ct(f) \cap D_r.      %%    \Ct(f)(r)
  \end{align*}

\claim For $t$ sufficiently small, $Y_t$
has exactly $2m$ points on the circle $S_r$.

To see this, consider the holomorphic function $\lambda_t(z) = h_t(r \cos z, r
\sin z)$.  The real zeros (mod $2\pi$) of $\lambda_t(z)$
correspond exactly to the zeros of $h_t(x,y)$ on the circle $S_r$.
For $\delta$ sufficiently small and $\alpha \in \R$ a non-root of
$\lambda_0$, the set $T = \{ z \st \Re(z) \in
[\alpha,\alpha+2\pi]\,, \Im(z) \in [-\delta,\delta]\}$ contains
exactly one representative (mod $2\pi$) of each of the $2m$ real roots of
$\lambda_0(z)$, and no other roots of $\lambda_0(z)$.  Observe
that the hypothesis that $\Ct(f)$ is not tangent to the circle
$S_r$ is precisely equivalent to the statement that the real zeros
of $\lambda_0(z)$ are all simple zeros.

For $t$ sufficiently small, by the argument principle
$\lambda_t(z)$ also has exactly $2m$ roots in $T$.  (Roots on the
boundary of $T$ are ruled out by the uniform convergence of
$\lambda_t \rightarrow \lambda_0$.) Hence $\lambda_t(z)$ has at
most $2m$ real roots (mod $2\pi$). However, again by the argument principle,
for any sufficiently small $\epsilon > 0$, if $t$ is sufficiently
small then the disk of radius $\epsilon$ around a real root of
$\lambda_0(z)$ contains exactly one root of $\lambda_t(z)$. Since
$\lambda_t(z)=0$ if and only if $\lambda_t(\overline{z})=0$, these
roots of $\lambda_t(z)$ must be real.  This establishes the claim;
in fact, we have shown that as $R \rightarrow \infty$, the point
$Q_j$ labeled $j$ on $Y_t$ approaches the point labeled $j$ on $Z$.
(See Figure~\ref{approaching-figure}.)

\smallskip

\begin{figure}[h]
\psfrag{Cf}{\Huge{\Red{$C_\theta(f)$}}}
\psfrag{Ch}{\Huge{\Blue{$C_\theta(h_t)$}}}
\psfrag{Sr}{\Huge{$S_r$}}
\psfrag{th}{\Huge{$\theta/n$}}
\psfrag{p0}{\Huge{\Red{$P_0$}}}
\psfrag{p1}{\Huge{\Red{$P_1$}}}
\psfrag{pn}{\Huge{\Red{$P_n$}}}
\psfrag{q0}{\Huge{\Blue{$Q_0$}}}
\psfrag{q1}{\Huge{\Blue{$Q_1$}}}
\psfrag{qn}{\Huge{\Blue{$Q_n$}}}
\resizebox{2.2in}{1.5in}{\includegraphics{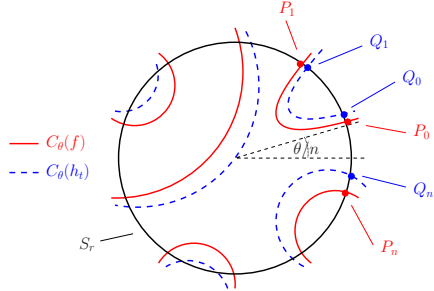}}
  \caption{How the curve $C_\theta(h_t)$ approaches the curve $C_\theta(f)$.
  \label{approaching-figure}}
\end{figure}

Let $C_1,\ldots,C_m$ denote the connected components of $\Ct(f)$, and let $\Delta > 0$ be
smaller than the distance between any two of the components: that is, if $z_j \in C_j$ and
$z_k \in C_k$ with $j \neq k$, then $|z_j-z_k| > \Delta$.  Take $R$ sufficiently large, as
above, so that $|P_j-Q_j| < \Delta/2$ for all $j$.

Let $\epsilon=\Delta/3$, and further take $t=1/R$ sufficiently small so as to
satisfy both parts of Lemma~\ref{lemma:technical}.  Note that the connected components of
$N_\epsilon(Z)$ are precisely $N_\epsilon(C_1)$, \dots, $N_\epsilon(C_m)$.

Let $U=\{z\st |z|<r\}$, and fix points $u_j \in C_j\cap U$ for $j=1,\dots,m$.
By part (2) of Lemma~\ref{lemma:technical}, there exist points $v_1,\dots,v_m\in
Y_t \cap U$ such that $|u_j-v_j|<\epsilon$.  That is,
$v_j\in N_\epsilon(C_j)$.  In particular, the $v_j$ are all distinct.

Let $C'_j$ be the connected component of $Y_t$ that
contains $v_j$; then $C'_j\subset N_\epsilon(Z)$ by part
(1) of Lemma~\ref{lemma:technical}. Since $C'_j$ is connected, we have $C'_j\subset
N_\epsilon(C_j)$.  In particular, the~$C'_j$ are all different.  Since
$Y_t$ has exactly $2m$ points on $|z|=r$, the~$C'_j$ exhaust
all of the connected components of $Y_t$. Finally, since
$|P_i - Q_i| < \Delta/2$, it must be the case that $P_i$ lies on~$C_j$
if and only if $Q_i$ lies on~$C'_j$. It follows that
$\Match(C_\theta(f),r) = \Match(C_\theta(g_R),r)$.
\end{proof}

We have now collected all the tools we need to prove Theorem~\ref{thm:matchings},
and thereby to characterize the combinatorial relationship between the matchings
$\Match(g_R,\theta)$ and $\Match(f,\theta)$.

\begin{proof}[Proof of Theorem \ref{thm:matchings}]
By Proposition~\ref{prop:disk}, for any $\delta > 0$, we may choose $r$ and $R$
sufficiently large that $\Match(\Ct(g_R),r)=\Match(f,\theta)$.  Further, take
$r$ and $R$ sufficiently large that
  \begin{itemize}
    \item $|\arg(f(z)/z^n)| < \delta$ for $|z| \ge r$,
    \item $\Match(\Ct(f),r) = \Match(f,\theta)$, and
    \item if $P_0,\ldots,P_{2n-1}$ and $Q_0,\ldots,Q_{2n-1}$ denote the points of
        $\Match(\Ct(f),r)$ and $\Match(\Ct(g_R),r)$ labeled $0,\ldots,2n-1$ respectively, then
    $$|\arg(P_k) - (k\pi + \theta)/n| < \frac{\delta}{2(n+1)} \ \ \text{and} \ \
      \left|\arg(Q_k)-\arg(P_k)\right| <\frac{\delta}{2(n+1)}\,.$$
  \end{itemize}

It follows from these properties that $\arg(g_R(z))$ differs from $n \arg(z) + \arg(z-R)$ by
at most $\delta$ if $|z| \ge r$ and $z \neq R$, and that $|\arg(Q_k) - (k\pi + \theta)/n| <
\delta/(n+1)$.

Recall that $\theta$ is assumed to lie in the interval $(0,\pi)$. We define the interval $I_k(\delta)$ for $0 \le k \le
2n$ as follows:
  \newcommand{\pad}{\rule[-6bp]{0bp}{6bp}}  %% to make the following case environment more legible
  $$I_k(\delta) = \begin{cases}
    (\frac{k\pi+\theta+\delta}{n+1},\ \frac{k\pi+\theta-\delta}{n})
      & \text{for} \ 0 \le k \le n-1, \pad\\
    (\frac{n\pi+\theta+\delta}{n+1},\ \pi)
      & \text{for} \ k=n, \pad\\
    (\frac{(k-1)\pi+\theta+\delta}{n},\ \frac{(k+1)\pi+\theta-\delta}{n+1})
      & \text{for} \ n+1 \le k \le 2n.
  \end{cases}$$
Observe that it is possible to take $\delta$ sufficiently small that all these intervals are
nonempty and $\delta<\theta$, and we do so.

\claim
If $\phi \in I_k(\delta)$ and $d \ge r$, then $z = de^{i\phi}$ does not lie on $\Ct(g_R)$.
That is, $\Ct(g_R)$ does not cross the half-line $\{de^{i\phi} \st d \ge r\}$.

To verify the claim, we must show that $\arg(g_R(de^{i\phi})) \not\equiv \theta \pmod{\pi}$.
Suppose first that $k \le n$, so that $\phi < \pi$.  Then $\arg(z-R) \in (\phi,\pi)$, and one
computes that
  $$n\arg(z) + \arg(z-R) \in (k\pi+\theta+\delta,\ (k+1)\pi + \theta - \delta) \,,$$
and therefore
  $$\arg(g_R(z)) \in (k\pi+\theta,\ (k+1)\pi+\theta)\,$$
as desired.
Similarly, if $k > n$, then $\phi > \pi$ and $\arg(z-R) \in (\pi,\phi)$, and we obtain the
same conclusion about $\arg(g_R(z))$.

\begin{figure}[h]
\psfrag{Cf}{\Huge{\Red{$C_\theta(f)$}}}
\psfrag{Ch}{\Huge{\Blue{$C_\theta(h_t)$}}}
\psfrag{Sr}{\Huge{$S_r$}}
\psfrag{R }{\Huge{$R\in\R$}}
\psfrag{p0}{\Huge{\Red{$P_0$}}}
\psfrag{p1}{\Huge{\Red{$P_1$}}}
\psfrag{pn}{\Huge{\Red{$P_{2n-1}$}}}
\psfrag{q0}{\Huge{\Blue{$Q_0$}}}
\psfrag{q1}{\Huge{\Blue{$Q_1$}}}
\psfrag{qn}{\Huge{\Blue{$Q_{2n-1}$}}}
\psfrag{I0}{\Huge{\DkGrn{$I_0$}}}
\psfrag{I1}{\Huge{\DkGrn{$I_1$}}}
\psfrag{I2}{\Huge{\DkGrn{$I_2$}}}
\psfrag{I2n}{\Huge{\DkGrn{$I_{2n}$}}}
\resizebox{4.0in}{4.0in}{\includegraphics{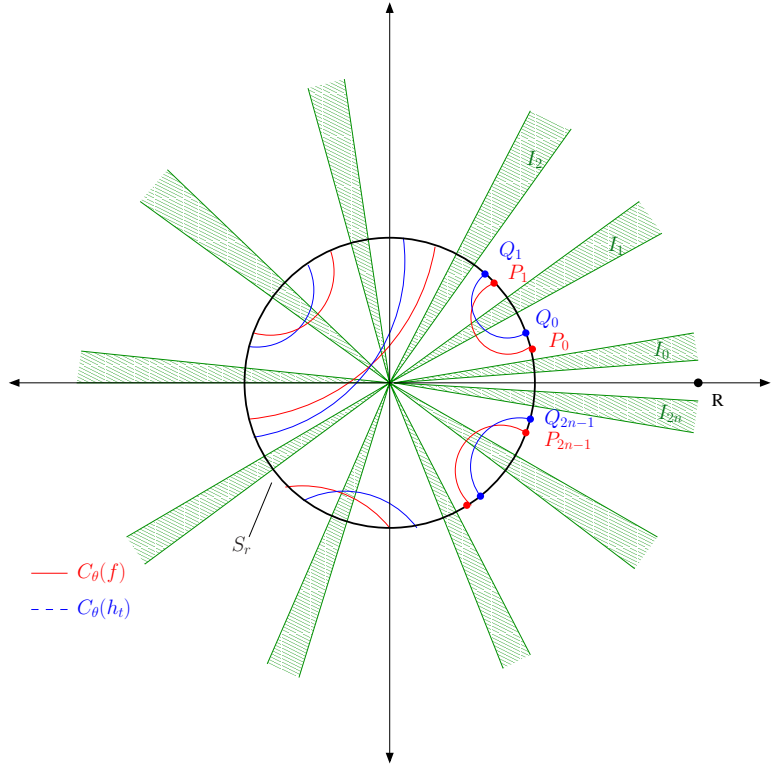}}
  \caption{Locating the new root $R$ between $I_{2n}$ and $I_0$
  \label{locatingR-figure}}
\end{figure}

\smallskip
Finally, we are in a position to complete the proof of Theorem \ref{thm:matchings}. Choose
$\phi_k \in I_k(\delta)$ for each $0 \le k \le 2n$.  One checks that
  $$(k\pi+\theta)/n, \ ((k+1)\pi+\theta)/(n+1) \in
    (\phi_k + \delta/(n+1) , \phi_{k+1} - \delta/(n+1))$$
for $0 \le k \le 2n-1$.  Moreover, $\theta/(n+1) < \phi_0$ and
$((2n+1)\pi+\theta)/(n+1) > \phi_{2n}$.

Since $|\arg(Q_k) - (k\pi + \theta)/n| < \delta/(n+1)$, we have established:
  \begin{itemize}
    \item $\arg(Q_k) \in (\phi_k,\phi_{k+1})$,
    \item the unique $j$ such that $(j\pi + \theta)/(n+1) \in (\phi_k,\phi_{k+1})$ is $j =
        k+1$, and
    \item $\Ct(g_R)$ does not cross the half-lines $\{de^{i\phi_k} \st d \ge r\}$ and
        $\{de^{i\phi_{k+1}} \st d \ge r\}$.
  \end{itemize}

Together, these facts imply that for $0 \le k \le 2n-1$, the connected component of $\Ct(g_R)$
containing $Q_k$ has an asymptote at angle $((k+1)\pi + \theta)/(n+1)$.  This proves that if
$a,b$ are matched in $\Match(f,\theta) = \Match(\Ct(g_R),r)$, then $a+1,b+1$ are matched in
$\Match(g_R,\theta)$.  By elimination, $0$ and $2n+1$ must also be matched in
$\Match(g_R,\theta)$; that is, $\Match(g_R,\theta) = \widehat{\Match(f,\theta)}$.
\end{proof}

We can now complete the proof of the Inverse Basketball Theorem.

\begin{proof}[Proof of Theorem~\ref{IBT}]
Suppose that $\alpha < \beta$.  Let $f_\eta(z)$ be the monic
polynomial $e^{in\eta} f(e^{-i\eta} z)$. Then the curve $C_{n\eta
+ \theta}(f_\eta)$ is the rotation of $\Ct(f)$ by $\eta$.
Taking $\eta$ such that $n\eta+\beta<\pi$, we see that
$\Basket(f,\alpha,\beta)=\Basket(f_\eta,0,\beta-\alpha)$; we can
therefore assume without loss of generality that $\alpha \neq 0$.

The set of all basketballs is closed under \defterm{rotation},
the operation of replacing each quartet $\{i,j\},\{k,\ell\}$
with $\{i-2,j-2\},\{k-2,\ell-2\}$ (with subtraction
taken modulo $4n$).  This transformation can also be realized on
the set of basketballs $\Basket(f,\alpha,\beta)$ arising from
monic polynomials.  Indeed, taking $\eta = -\pi/n$, we see that
$a,b$ are matched in $\Match(f_{\pi/n},\theta)$ if and only if
$a-1,b-1$ are matched in $\Match(f,\theta)$.  Since this is
true for both $\theta=\alpha,\beta$, we see that
$\Basket(f_{\pi/n},\alpha,\beta)$ is the rotation of
$\Basket(f,\alpha,\beta)$.

The basketballs are also closed under \textit{half-rotation},
the operation that replaces each quartet $\{i,j\},\{k,\ell\}$ with
$\{k-1,\ell-1\},\{i-1,j-1\}$.  Choose $0 < \gamma < \beta-\alpha$,
and take $\eta = -(\alpha+\gamma)/n$.  One checks similarly that
the basketball $\Basket(f_\eta,\beta-\alpha-\gamma,\pi-\gamma)$ is
the half-rotation of $\Basket(f,\alpha,\beta)$.

We proceed by induction on the degree of $f$, the case $n=1$ being
trivial. Let $B$ be any basketball of order $n+1$.  This contains
an ear by Proposition~\ref{thm:ear-thm}, and by rotation we may
assume that the ear contains the pair $\{4n+3,1\}$.  The other
pair in this ear is either $\{4n+2,0\}$ or $\{0,2\}$.  If the ear
is $\{4n+3,1\},\{0,2\}$ then the half-rotation of $B$ contains the
ear $\{4n+3,1\},\{4n+2,0\}$.  By the observation in the previous
paragraph, we may assume without loss of generality that the ear is
$\{4n+3,1\},\{4n+2,0\}$.

Let $(M_{\alpha},M_{\beta})$ be the ordered pair of matchings on $[0,2n+1]$
corresponding to $B$.  By assumption, the pair $\{0,2n+1\}$ is
contained in both $M_\alpha$ and $M_\beta$.

For $\theta = \alpha,\beta$, let $M'_{\theta}$ be the matching on
$[0,2n-1]$ such that $\{a,b\}\in M'_{\theta}$ if and
only if $\{a+1,b+1\}\in M_{\theta}$. Then the
bimatching $B'$ corresponding to $(M'_{\alpha},M'_{\beta})$ is a
basketball of degree $n$, and so by the induction hypothesis
$B'=\Basket(f,\alpha,\beta)$ for some monic polynomial $f(z)$.
Since $\alpha \neq 0$, by Theorem \ref{thm:matchings}
we see that $B$ is the basketball of
$(z-R)f(z)$ for $R$ sufficiently large, and we are done.
\end{proof}

%===========================================================================

\section{Necklaces of matchings} \label{necklace-section}

More generally, we are
interested in classifying the possibilities for the topology of the family
  $$\CC(f)=\{C_\theta(f) \st \theta\in \R/\pi\Z\},$$
which is truly an invariant of $f$ itself, not depending on a choice of angle.
The family $\CC(f)$ is fibered over the base $\R/\pi \Z$ with fiber $C_{\theta}(f)$
above $\theta$.  In particular, suppose that $f(z)$ has distinct roots.
Let $z_1,\ldots,z_{n-1}$ be the roots of $f'$, and suppose further
that $\arg f(z_1), \ldots, \arg f(z_{n-1})$ are distinct.
Then the family $\CC(f)$ has $n-1$ fibers with ordinary singularities,
and is smooth elsewhere; that is, it is smooth over $n-1$ open arcs
$\mathcal{A}_i$ arranged cyclically around the circle $\R/\pi\Z$.
The noncrossing matching $M(f,\theta)$ is the same for all $\theta \in
\mathcal{A}_i$; denote it by $M_i$.  We thus obtain an $(n-1)$-tuple
$\MM=(M_1,\ldots,M_{n-1})$ of noncrossing matchings, and the data of $\MM$
determines the topology of $\CC(f)$.

Let $M_n$ be the matching $\left\{\{i-1,j-1\} \st \{i,j\} \in
M_1\right\}$, taking all indices modulo $2n$ as usual.  Then the
$(n-1)$-tuple of matchings $(M_1,\dots,M_{n-1})$
has the following property: for $0 \le t
\le n-1$, the matching $M_{t+1}$ can be obtained from $M_t$ by
taking a suitable pair of pairs $\{i,j\},\{k,\ell\}$ and replacing
them with $\{i,\ell\},\{j,k\}$.  Call an $(n-1)$-tuple of
noncrossing matchings possessing this property a \emph{necklace of
matchings of order~$n$}.  For example, from this point of view the
basketball of a
quadratic polynomial $f$ is determined by the two possibilities for
the necklace $(M_1)$, and by whether $\arg f(z_1) \pmod{\pi}$ is
greater or smaller than $\pi/2$ (see also Remark~\ref{rmk:discriminant}).

It is natural to ask whether every necklace of order $n$ must arise from a
polynomial of degree $n$.  To that end, define a \emph{multiear} to be an
integer $i$ for which
$\{i,i+1\}\in M_1,\ldots,M_t$ and $\{i-1,i\} \in M_{t+1},\ldots,M_n$
for some $t$.  If every necklace contains a multiear, then Theorem
\ref{thm:matchings} can be used to show that every necklace arises from a
polynomial; otherwise, new techniques will be necessary.  We note the following
result, obtained via exhaustive computer calculation:

\begin{prop} For $n \le 8$, we have:
\begin{enumerate}
\item Every necklace of order $n$ contains a multiear;

\item If $\MM=(M_1,\ldots,M_{n-1})$ is a necklace of order $n$ and $M_n$ is
defined as above, then for all $1 \le t < u \le n$ the bimatching
corresponding to $(M_t,M_u)$ is a basketball;

\item The number of necklaces of order $n$ is $2(2n)^{n-2}$.
\end{enumerate}
\end{prop}

We remark that (1) implies (2).  Indeed, (1) implies that every necklace of
order at most $8$ arises from a polynomial, and (2) is automatically satisfied
by any necklace arising from a polynomial.  The numbers $2(2n)^{n-2}$ appear
as sequence A097629 in \cite{EIS}, enumerating unrooted directed trees on $n$ vertices;
we do not know a bijective reason for this.

%===========================================================================

\section*{Acknowledgments}

It is our pleasure to thank Mira Bernstein, Pete Clark, Valery Liskovets,
Grisha Mikhalkin, Vic Reiner and Timothy Walsh for helpful conversations.
We began our work on this project while the first two authors were instructors
and the third author was a student at Canada/USA Mathcamp, a
summer program for high school students.  During much of the time
that this work was ongoing, the first two authors were
postdoctoral fellows at the University of Minnesota and McGill
University, respectively.

\end{document}